\documentclass[11pt,twoside]{article}
\usepackage{latexsym}
\usepackage{amssymb,amsbsy,amsmath,amsfonts,amssymb,amscd}

\usepackage{epsfig, graphicx}
\usepackage{epstopdf}
\usepackage{graphicx}
\usepackage{subfigure}
\usepackage{epstopdf}
\usepackage{xcolor}
\newcommand{\commentout}[1]{}

\newcommand {\Chi} {{\bf \raise 2pt \hbox{$\chi$}} }
\newcommand {\f}   {\frac}
\newcommand {\p}   {\partial}

\newcommand {\proof} {\noindent {\bf Proof}. }
\newcommand{\beq}{\begin{equation}}
\newcommand{\eeq}{\end{equation}}
\newcommand{\bea} {\begin{array}{rl}}
\newcommand{\eea} {\end{array}}
\newcommand{\bepa}{\left\{ \begin{array}{l}}
\newcommand{\eepa} {\end{array}\right.}
\newtheorem{theorem}{Theorem}[section]
\newtheorem{lemma}[theorem]{Lemma}

\newtheorem{remark}[theorem]{Remark}

\newcommand{\qed}{{ \hfill
                       {\unskip\kern 6pt\penalty 500 \raise -2pt\hbox{\vrule\vbox to 6pt{\hrule width 6pt
                       \vfill\hrule}\vrule} \par}   }}
\title{\Large \bf  Uniform convergent scheme for strongly anisotropic diffusion
equations with closed field lines}

\author{Yihong Wang
\thanks{Institute of natural sciences and department of mathematics,
Shanghai Jiao Tong University, Shanghai, 200240, China.}
\and Wenjun Ying
\thanks{Institute of natural sciences and department of mathematics,
Shanghai Jiao Tong University, Shanghai, 200240, China.}
\and Min Tang
\thanks{Institute of natural sciences and department of mathematics,
Shanghai Jiao Tong University, Shanghai, 200240, China.}}

\date{\today}

\begin{document}
\maketitle
\pagestyle{plain}
\pagenumbering{arabic}
\begin{abstract}
In magnetized plasma, the magnetic field confines particles around field lines. The ratio between the intensity of the parallel and perpendicular viscosity or heat conduction may reach the order of $10^{12}$. When the magnetic fields have closed field lines and form a ``magnetic island'', the convergence order of most known schemes depends on the anisotropy strength. In this paper, by integration of the original differential equation along each closed field line, we introduce a simple but very efficient asymptotic preserving
reformulation, which yields uniform convergence with respect to the anisotropy strength. Only slight modification to the original code is required and neither change of coordinates nor mesh adaptation is needed. Numerical examples demonstrating the performance of the new scheme are presented.

 \end{abstract}
\\[3mm]


\section{Introduction}\label{intro}

Anisotropic diffusion is encountered in many fields such as heat conduction in magnetized plasma \cite{Bram 14}, flows in porous media \cite{Ashby 99}, image processing or oceanic flows \cite{Galperin 10,Treguier 02}. We are particularly interested in the anisotropic diffusion in magnetized plasmas. It has extremely anisotropic diffusion tensors of heat conduction in fusion plasmas, where the particles are confined by the magnetic field and the particle diffusion is much faster along the field lines than in the perpendicular direction \cite{Bram 14}.

In magnetized plasma, magnetic field lines can be open or closed. 
The region that has closed field lines inside is called a {``}magnetic island". The field lines outside the region are open. Plasma transportation across these closed field lines depends on the perpendicular diffusion coefficients. Small diffusion coefficients in the perpendicular direction lead to plasma confinement. The closed field lines relate to the important physical process {``}magnetic reconnection" and can appear as a consequence of instabilities or external perturbations.

The problem under consideration is the following two-dimensional diffusion equation with anisotropic diffusivity:
\beq\label{eq:ellipT}
\left\{\begin{array}{ll}
-\nabla\cdot (A\nabla u^{\epsilon})=f,&\mbox{in $\Omega$},\\
u^{\epsilon} = g  ,&\mbox{on $\partial\Omega$},
\end{array}
\right.
\eeq
where $\Omega\subset R^2$ is a bounded domain with boundary $\p \Omega$. The
diffusion tensor is given by
 \beq\label{eq:A}
A(x,y)=\left(\begin{array}{cc}\cos\theta&-\sin\theta\\ \sin\theta&\cos\theta\end{array}\right)
\left(\begin{array}{cc}1/\epsilon&0\\0&\alpha\end{array}\right)\left(\begin{array}{cc}\cos\theta&\sin\theta\\-\sin\theta&\cos\theta\end{array}\right),
\eeq where $\alpha$ is $O(1)$ and the parameter $0<\epsilon<1$ can be very small. All coefficients $\alpha$, $\epsilon$ and $\theta$ may depend on space.
The direction of the anisotropy or the magnetic field line is given by a unit vector field $\mathbf{b}=(\cos\theta,\sin\theta)^T$, while $\alpha$ and $1/\epsilon$ are respectively the perpendicular and parallel diffusion coefficients.
The problem becomes highly anisotropic when $\epsilon\ll1$.

Numerical simulations for anisotropic diffusion problems have been addressed by a lot of researchers and engineers; see the review in \cite{Herbin 08}.  Methods used today include finite volume method \cite{Potier 05,Sheng 09, Yuan 08}, finite difference method \cite{Bram 14}, mimetic finite difference method \cite{Hyman 02}, discontinuous Galerkin method \cite{ARNOLD82,ALEXANDRE 09},
finite element method \cite{Gunter07, Hou 97, Li 10, Pasdunkorale 05} and so on. These methods are usually efficient for a selected range of $\epsilon$ but loss convergence when $\epsilon\ll h$ ($h$ is the mesh size).

 A field-aligned coordinate system is usually employed for plasma simulations. However, it may run into problems when there are magnetic re-connections or
highly fluctuating field directions. Schemes with non-aligned meshes in case of varying anisotropy
have been studied as well, for instance \cite{Bram 14,Gunter07,Gunter05} and the references there in. Since the numerical errors in the direction parallel to the magnetic field may have significant effect on the perpendicular direction,
several difficulties arise for strongly anisotropic diffusion problems with non-aligned meshes, which include pollution on the perpendicular direction by the parallel diffusion and loss of convergence, etc \cite{Bram 14}.

In this paper, we are interested in the case when the magnetic field exhibits closed field lines.
First, due to the closed field lines, numerical discretizations of the original problem using magnetic field aligned coordinates lead to very badly conditioned systems when $\epsilon$ becomes small. The limiting model is not well-posed and admits infinitely many solutions \cite{Narski13}, as adding any function constant in the region covered by closed field lines keeps a solution to the limiting model. Second, when non-aligned coordinates are used, most known schemes loss convergence when $\epsilon \to 0$ \cite{Bram 14,Gunter07,Gunter05}. If the field line directions are taken into account in the numerical discretization, uniform convergence can be achieved when the closed field line is symmetric, but all schemes in \cite{Bram 14,Gunter07,Gunter05} fail to converge when $\epsilon\ll h$ in the tilted closed field line case.

That the scheme convergence is independent of $\epsilon$ can be considered as an Asymptotic preserving (AP) property.
Asymptotic preserving (AP) methods for strongly anisotropic diffusion have been studied in a series paper by Degond et.al  \cite{Degond101,Degond102,Degond121,Degond122,Negulescu14,Mentrelli12}. AP schemes have the advantage that the condition number does not scale with the anisotropy, when using Neumann or periodic boundary conditions, thus can deal with all $\epsilon$ ranging from $O(1)$ to very small. The main idea in \cite{Degond101,Degond102,Degond121, Degond122, Negulescu14, Mentrelli12} are based on macro-micro decomposition, and reformulate the original equation into a system that keeps well-posed when $\epsilon\to 0$. For the closed field line case, Narski and Ottaviani \cite{Narski13} developed a uniform convergent scheme by introducing a penalty stabilization term, where a tuning parameter is needed.


In this paper, we propose a different approach, which is simple and easily extendable.
The main idea is that we cut each closed field line at some point $(x_{0}, y_{0})$, which
is treated as the begin and end points of the cut field line.
Then instead of discretizing the equation at the point  $(x_{0}, y_{0})$ locally, we work with integration of the  differential equation along the cut field line. By using continuity conditions at the cutting point $(x_0,y_0)$, the singular $1/\epsilon$ terms disappear, so that the reformulated problem becomes well-posed.

The idea is inspired by the AP method developed in \cite{tang 16} for the strongly anisotropic diffusion equation with Neumann boundary conditions.
After replacing one of the Neumann boundary condition by integration of the original differential equation along the
 field line, an equivalent new system that is well-posed in the limit $\epsilon\to 0$ can be constructed. Any numerical discretization based on this new system has uniform convergence with respect to $\epsilon$.

 The main advantage of the field line integration approach is that: 1) uniform convergence is numerically achieved even for the tilted closed field line case; 2) only slight modification to the original code is needed, which make it attractable for engineers; 3) no tuning parameter is introduced; 4) no coordinate change or mesh adaptation is required.

%

The paper is organized as follows. In section 2, we illustrate the idea of the reformulation. Numerical discretizations are given in Section 3, where
nine point finite difference method (FDM) is used for the reformulated system.
 Several numerical examples are presented in section 4, which show the uniform convergence of the scheme with respect to the anisotropy. Finally we conclude with some discussions in section 5.

\section{The Asymptotic Preserving Reformulation}

\begin{figure}[htb]
\centering
{
\includegraphics[width=7.5cm] {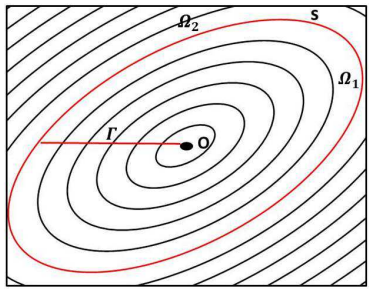}}
\caption{
Closed
field lines.}\label{fig:branch close}
\end{figure}

Let $\Omega = (-a, a) \times (-b, b)$ be the computational domain. Assume it can be partitioned into two non-overlapping open subdomains $\Omega_1$ and $\Omega_2$ by a separatrix $S$ such that
$\Omega = \Omega_1 \cup \Omega_2\cup S$. All field lines in $\Omega_1$ are closed and those in $\Omega_2$ are open.

We first discuss about the subdomain $\Omega_2$ with open field lines.
In the open field line subdomain $\Omega_2$, Eqn \eqref{eq:ellipT} remains the same, any scheme in \cite{Bram 14,Gunter07,Gunter05} can be applied.
Whether a field line $\ell$ is considered as an open or closed field line depends on whether it intersects with the boundary $\partial\Omega$ or not. If $\partial \Omega\cap \ell =\emptyset$, then $\ell$ is a closed field line. If 
$\partial \Omega\cap \ell\neq\emptyset$, then $\ell$ is open. The region inside $S$ is a "magnetic island". Outside $S$, the field lines are open (see Figure \ref{fig:branch close} for an illustration).


In the subsequent part, we will construct an equivalent system in $\Omega_1$ ($\Omega_1\cup S$).
To illustrate the idea, we use the notations in Figure \ref{fig:branch close}.  It is important to note that when there  are closed field lines, one can find a singular point inside $\Omega_1$  where $\cos\theta$ and $\sin\theta$ have no definition. It  is considered as the center of the "magnetic island". Suppose $O(0,0)$ is the center of the "magnetic island" and $\ell$ $(l\in \Omega_1\cup S)$ is a close field line.

\subsection{An interface problem}

Let $\mathbf{b}=(\cos\theta,\sin\theta)^T$, $\mathbf{b}_\bot=(-\sin\theta,\cos\theta)^T$, Eqn \eqref{eq:ellipT} can be written as
\beq\label{eq:SimellipA}
- ( \mathbf b \cdot \nabla) (   \frac{1}{\epsilon} \mathbf   b^{} \cdot \nabla u^{\epsilon}) -
 ( \nabla \cdot\mathbf b) (   \frac{1}{\epsilon}\mathbf b^{} \cdot \nabla u^{\epsilon}) - \nabla
   \cdot \big( \alpha \mathbf b_{\perp} (\mathbf b_{\perp} \cdot\nabla u^{\epsilon})\big) = f.
   \eeq
   We can pick any $\mathbf{b}$-field line $\ell$ and parameterize it by the arc length $s$.  Accordingly, $\partial_s$ will denote the derivative
along the line $\ell$, i.e. $\partial_s = \mathbf{b} \cdot \nabla$. Then, {Eqn} \eqref{eq:SimellipA}
can be written as
\beq\label{eq:rew1}
 -  \partial_s(\frac{1}{\epsilon}\p_s u^{\epsilon}) -  ( \nabla
   \cdot \mathbf{b}) \frac{1}{\epsilon} \partial_s^{} u^{\epsilon} -  \nabla
   \cdot \big( \alpha \mathbf{b}_{\perp} (\mathbf{b}_{\perp}\cdot \nabla u^{\epsilon})\big)  = f.
   \eeq
In the strong anisotropic diffusion limit as $\epsilon\rightarrow 0$, \eqref{eq:rew1} yields
\beq\label{eq:rew2lim}
- \partial_s^2 u^0 -  ( \nabla
   \cdot \mathbf{b}) \partial_s^{} u^0 =0.
\eeq
If $\ell$ is closed, any function that is constant along $\ell$ satisfies \eqref{eq:rew2lim}. There is no uniqueness of the solution to the limiting model, which explains the ill-posedness when a field-aligned coordinate is used in the numerical discretization.

We cut $\ell$ at a point $(x_{0},y_{0})$. Starting from $(x_0,y_0)$ and following the direction $\mathbf{b}$, we can determine $(x_{0^+},y_{0^+})$ and $(x_{0^-},y_{0^-})$, which are respectively the start and end point of $\ell$. Since $\ell$ is closed, \beq(x_{0^+},y_{0^+})=(x_{0^-},y_{0^-})=(x_{0},y_{0}).\label{eq:x0y0}\eeq

When $A(x,y)$ is continuous, the solution to \eqref{eq:ellipT} belongs to $C^1(\Omega)$, thus both $u^\epsilon$ and $\nabla u^\epsilon$ are continuous at $(x_{0},y_{0})$. However, this regularity requirement is not necessarily true for the limiting equation \eqref{eq:rew2lim}. Since $s$ is the arc length,
$s = 0$ corresponds to $(x_{0^+},y_{0^+})$ and $s = L_{\ell}$ with
$L_{\ell}$ being the length of $\ell$ corresponds to
$(x_{0^-},y_{0^-})$.  {Eqn} \eqref{eq:rew2lim} holds at all points on $\ell$, which indicates that
$$u^0|_{s=0}=u^0|_{s=L_{\ell}}, \qquad
\partial_s^{} u^0 |_{s=0}=\partial_s^{} u^0 |_{s=L_{\ell}}.$$
Therefore, due to the regularity requirement, the following connection conditions
\begin{equation}\label{eq:connection}
u^{\epsilon}_+=u^{\epsilon}_-,\qquad
\frac{1}{\epsilon} \mathbf{t}\cdot\mathbf{b}  (\mathbf{b}\cdot \nabla u_{+}^{\epsilon}) + \alpha \mathbf{t}\cdot\mathbf{b}_{\perp}
   (\mathbf{b}_{\perp}\cdot \nabla u_{+}^{\epsilon}) =  \frac{1}{\epsilon} \mathbf{t}\cdot\mathbf{b}  (\mathbf{b}\cdot \nabla u_{-}^{\epsilon}) + \alpha \mathbf{t}\cdot\mathbf{b}_{\perp}
   (\mathbf{b}_{\perp}\cdot \nabla u_{-}^{\epsilon}),\end{equation}
should be satisfied independent of $\epsilon$. Here
$u^{\epsilon}_+=u^{\epsilon}(x_{0^+},y_{0^+})$, $u^{\epsilon}_-=u^{\epsilon}(x_{0^-},y_{0^-})$ and
\eqref{eq:connection} can be written as
$$u^{\epsilon}(x_{0^+},y_{0^+})=u^{\epsilon}(x_{0^-},y_{0^-}),\qquad \mathbf{t}\cdot A(x_0,y_0)\nabla u^{\epsilon}(x_{0^+},y_{0^+})=\mathbf{t}\cdot A(x_0,y_0)\nabla u^{\epsilon}(x_{0^-},y_{0^-}),
$$
where $\mathbf{t}$ is the tangent direction of the field line across $(x_0,y_0)$.

On each closed field line, we choose a cutting point. To describe the calculation better, we can suppose all cutting points are on the negative X axis. Let $\Gamma$ be the set of all cutting points (see Figure \ref{fig:branch close}). Then we have the following  interface problem:
\begin{equation}\label{eq:interface}
\left\{\begin{aligned}
&- ( \mathbf b \cdot \nabla) (   \frac{1}{\epsilon} \mathbf   b^{} \cdot \nabla u^{\epsilon}) -
 ( \nabla \cdot\mathbf b) (   \frac{1}{\epsilon}\mathbf b^{} \cdot \nabla u^{\epsilon}) - \nabla
   \cdot \big( \alpha \mathbf b_{\perp} (\mathbf b_{\perp} \cdot\nabla u^{\epsilon})\big) = f,\quad &\mbox{on $\Omega_1\setminus\Gamma$},\\
   &u^{\epsilon}(x_{0^+},y_{0^+})=u^{\epsilon}(x_{0^-},y_{0^-}),\quad &  \mbox{for $\forall (x_0,y_0)\in\Gamma$}, \\
  & \mathbf{t}\cdot A(x_0,y_0)\nabla u^{\epsilon}(x_{0^+},y_{0^+})=\mathbf{t}\cdot A(x_0,y_0)\nabla u^{\epsilon}(x_{0^-},y_{0^-}), \quad &  \mbox{for $\forall (x_0,y_0)\in\Gamma$}.
\end{aligned}
\right.	
\end{equation}
\begin{remark}
If $u^\epsilon$ is the solution to the original problem \eqref{eq:ellipT}, it satisfies the interface problem \eqref{eq:interface}. On the other hand, it is easy to check that when $u^\epsilon\in C^1(\Omega)$ and satisfies the interface problem \eqref{eq:interface}, it is a solution to the original problem \eqref{eq:ellipT}.	
\end{remark}

\subsection{Reformulation of the interface problem}

By similar argument as in section 2.1, the interface problem \eqref{eq:interface} is ill-posed when $\epsilon\to 0$. However, we can replace the continuity of the derivatives at the interface points by integration of the equation along the field lines, and derive an equivalent system that is well-posed when $\epsilon\to 0$. The details are as follows.

First of all, we introduce a function $E$ defined on {$\Omega_2$}. Assume the field line $\ell\subset \Omega_2$ starts and ends at $(x_0,y_0)$. The function $E$ solves on the field line $\ell$ the differential equation
\beq\label{eq:E}  \partial_s E = ( \nabla \cdot \mathbf{b}) E,\qquad E(0) = 1.   \eeq
The solution to \eqref{eq:E} is
${E}(s)=e^{\int_0^ {s} \nabla \cdot \mathbf{b} \,\,ds'}$.
The following lemma indicates that if the limiting solution $u^0$ exists, $E(L_{\ell}) = E(0) = 1$. Here, $L_{\ell}$ is the arc length
of the field line $\ell$. 

\begin{lemma}
If there exists a function $v(x,y)\in C^2(\ell)$ such that $\mathbf{b}=\frac{1}{| \nabla v|}(-v_y, v_x)^T$, $\int_0^{L_{\ell}}\nabla\cdot \mathbf {b}\, {ds} =0$.
\end{lemma}
\proof: Let $\mathbf{b}=(\dot{x}(t),\dot{y}(t))^T$
with
\begin{displaymath}
\dot{x}(t)=\f{-v_y(x(t), y(t))}{| \nabla v|},  
\quad 
\dot{y}(t)=\f{v_x(x(t), y(t))}{| \nabla v|}, 
\end{displaymath}
$t=t_0$ correspond to $s=0$ and $t=t_l$ correspond to $s=L_{\ell}$. 
Therefore, $t=t_0$ and $t=t_\ell$ correspond to the same point in $\Omega$ and $v_x\mid_{t=t_0}=v_x\mid_{t=t_l}$, $v_y\mid_{t=t_0}=v_y\mid_{t=t_l}$. Besides from $\dot{x}(t)^2+\dot{y}(t)^2=1$, $ds=dt$. Then
 \beq
 \begin{aligned}
\int_0^{L_{\ell}}\nabla \cdot \mathbf {b}^{}\,{ds}
&=\int_{t_0}^{t_l}\partial_x
\biggl (\f{-v_y}{| \nabla v|} \biggr )+\partial_y \biggl (\f{v_x}{| \nabla v|} \biggr ){dt}
\\[4pt] 
\nonumber 
&=\int_{t_0}^{t_l}\f{-v_x^2 v_{xy}-v_x v_y v_{yy}+v_y^2v_{xy}+v_xv_y v_{xx}}{(v_x^2+v_y^2)^{3/2}}{dt}\\
&=\int_{t_0}^{t_l}-\f{(v_y v_{xy}+v_xv_{xx})\dot{x}(t)+(v_x v_{xy}+v_y v_{yy})\dot{y}(t)}{v_x^2+v_y^2}{dt}\\
&=\int_{t_0}^{t_l}-\f{\f{1}{2}(v_x^2+v_y^2)_x\dot{x}(t)+\f{1}{2}(v_x^2+v_y^2)_y\dot{y}(t)}{v_x^2+v_y^2}{dt}.
\\
&=\int_{t_0}^{t_l}-\f{d}{dt} \biggl \{\f{\ln(v_x^2+v_y^2)}{2} \biggr \}dt=0.
\end{aligned}
\eeq
\qed

If the limiting solution $u^0$ exists, from \eqref{eq:rew2lim}, it satisfies $\mathbf {b}\cdot \nabla u^0=0$ along the closed field lines. Then the vector field $\mathbf {b}$ can be determined by $ \mathbf{b}=\frac{1}{|\nabla u^0|}(-u^0_y, u^0_x)^T$. The above Lemma gives $E(0)=E(L_{\ell})=1$. It is important to note that this property does not rely on the value of $\epsilon$.

 We multiply both sides of \eqref{eq:rew1} by $E$ and combine the two singular $O(1/\varepsilon)$ terms to get 
\beq\label{eq:rew2}
-  \partial_s^{} ( E \frac{1}{\epsilon}\partial_s^{} u^{\epsilon}) - E \nabla
   \cdot \big( \alpha \mathbf{b}_{\perp} (\mathbf{b}_{\perp}\cdot \nabla u^{\epsilon})\big)  = {Ef}.
   \eeq
 Then integrating \eqref{eq:rew2} over $\ell$ gives us 
 \beq\label{eq:rew31}
 -E\frac{1}{\epsilon}\partial_s^{} u^{\epsilon}\big|_{s = 0}^{s = L_{\ell}} - \int_0^{L_{\ell}} E \nabla
   \cdot \big( \alpha \mathbf{b}_{\perp} (\mathbf{b}_{\perp}\cdot \nabla u^{\epsilon})\big)  {ds} = \int_0^{L_{\ell}} Ef{ds}.
\eeq
From $\partial_s^{} u^\epsilon |_{s=0}=\partial_s^{} u^\epsilon |_{s=L_{\ell}}$, and $E(0)=E(L_{\ell})=1$, let 
$S_\Gamma$ be the set of all closed field lines that satisfy $\ell(L_\ell)\in\Gamma$,
\eqref{eq:interface} can then be reformulated 
 as
  \beq\label{eq:SimellipC}
\left\{\begin{aligned}
&- ( \mathbf b \cdot \nabla) (   \frac{1}{\epsilon} \mathbf   b^{} \cdot \nabla u^{\epsilon}) -
 ( \nabla \cdot\mathbf b) (   \frac{1}{\epsilon}\mathbf b^{} \cdot \nabla u^{\epsilon}) - \nabla
   \cdot \big( \alpha \mathbf b_{\perp} (\mathbf b_{\perp} \cdot\nabla u^{\epsilon})\big) = f,\quad &\mbox{on $\Omega_1$},\\
   &u^{\epsilon}(x_{0^+},y_{0^+})=u^{\epsilon}(x_{0^-},y_{0^-}),\quad &  \mbox{for $\forall (x_0,y_0)\in\Gamma$}, \\
  &\int_0^{L_{\ell}} E \Big(\nabla
   \cdot \big( \alpha \mathbf{b}_{\perp} (\mathbf{b}_{\perp}\cdot \nabla u^{\epsilon})\big)+f\Big)  {ds} = 0\quad &  \mbox{for $\forall \ell\in S_\Gamma$}.
\end{aligned}
\right.	
\eeq
It is easy to check that when $\epsilon$ is finite, the above system is equivalent to the interface problem \eqref{eq:interface}.
The advantage of the reformulation \eqref{eq:SimellipC} is when $\epsilon\to 0$, its leading order gives
  \beq\label{eq:SimellipD}
\left\{\begin{aligned}
&  - \partial_s^2 u^0 - ( \nabla
   \cdot \mathbf{b}) \partial_s^{} u^0 = 0,& \mbox{in $\Omega_1$},\\
&  u^0\mid_{s=0}=u^0\mid_{s=L_{\ell}},  &\mbox{on $\Gamma$},\\
&\int_0^{L_{\ell}} E \Big(\nabla
   \cdot \big( \alpha \mathbf{b}_{\perp} (\mathbf{b}_{\perp}\cdot \nabla u^{0})\big)+f\Big)  {ds} = 0,  &\quad \mbox{for $\forall \ell\in S_\Gamma$}.
\end{aligned}
\right.
\eeq
The first two equations in \eqref{eq:SimellipD} indicate that $u^0$ is constant along $\ell$ while the third equation determines the constant.



\bigskip
In summary, combing the equation on $\Omega_2$ with \eqref{eq:SimellipC}, we get the asymptotic preserving reformulation for \eqref{eq:ellipT} on the whole computational domain such that
  \beq\label{eq:ALLt}
\left\{\begin{aligned}
&- ( \mathbf b \cdot \nabla) (   \frac{1}{\epsilon} \mathbf   b^{} \cdot \nabla u^{\epsilon}) -
 ( \nabla \cdot\mathbf b) (   \frac{1}{\epsilon}\mathbf b^{} \cdot \nabla u^{\epsilon}) - \nabla
   \cdot \big( \alpha \mathbf b_{\perp} (\mathbf b_{\perp} \cdot\nabla u^{\epsilon})\big) = f,\quad &\mbox{on $\Omega\setminus\Gamma$}, \\
   &u^{\epsilon}(x_{0^+},y_{0^+})=u^{\epsilon}(x_{0^-},y_{0^-}),\quad &  \mbox{on $\Gamma$},  \\
  & \int_0^{L_{\ell}} E \Big(\nabla
   \cdot \big( \alpha \mathbf{b}_{\perp} (\mathbf{b}_{\perp}\cdot \nabla u^{\epsilon})\big)+f\Big)  {ds} = 0,\quad &  \mbox{on $\ell\in S_\Gamma$}, \\
&u^\epsilon=g,\quad &\mbox{on $\p \Omega$},
\end{aligned}
\right.
\eeq
with $E$ defined as in \eqref{eq:E}.
\begin{remark}
Although $\mathbf{b}$ is continuous, {$\nabla\cdot\mathbf{b}$} may have jumps. This is always the case for ``magnetic islands". Examples can be found in our numerical examples in section 4. Then, $E$ are only $C^0$ but not $C^1$. To get a good approximation of the integration in \eqref{eq:ALLt}, we have to take into account the discontinuities in {$\nabla\cdot\mathbf{b}$} .	
\end{remark}

\section{Numerical Discretization}

Suppose the rectangle domain $\Omega=(-a,a)\times(-b,b)$ is partitioned into a 
uniform Cartesian grid with nodes
 $$
\begin{array}{ll}
\mathbf{z}_{i,j}=(x_i,y_j),\qquad&\mbox{ $i=0,\pm 1,\cdots \pm I$;
	$j=0,\pm 1,\cdots \pm J$}.
\end{array}
$$
Here, $I$ and $J$ are two positive integers,
$x_{i}= ih_x$ and $y_{j}= jh_y$ with  $h_x=a/I$, $h_y=b/J$.

The major difficulty in discretizing the system \eqref{eq:ALLt} lies in the 
integration of the differential equations along closed field lines. 
First we explain details on numerical calculation of closed field lines and 
selection of the quadrature points for numerical integration. Then we present 
the scheme that is used in our numerical tests. However, the reformulation 
does not rely on the specific discretization. Other schemes can be applied as well.

\subsection{Determine the closed field lines and the quadrature}
\textbf{Determine the open field lines: }
Based on the vector field $\mathbf{b}=(\cos\theta,\sin\theta)^T$, we can determine the field line by the following nonlinear ODE system:
 \beq\label{eq:ODE}
\left\{\begin{array}{ll}
\displaystyle \dot x(t)=\cos\theta(x(t),y(t)),\quad x(0)=x_0,
\\
\displaystyle \dot y(t)=\sin\theta(x(t),y(t)),\quad y(0)=y_0.
\end{array}
\right.
\eeq
To get the field line $\ell$,
we use a high order Runge-Kutta method with small time step to solve the above system. Numerically, we can successively obtain the discrete points
 $\{ (x_0,,y_0), \, (x_1, y_1), \, \cdots, (x_n, y_n) \, \}$ which form a discretization of the corresponding field line. If the field line is open, the calculations stop when $x_n$ or $y_n$ exceeds the boundary.
\begin{figure}[htb]
\centering
{
\subfigure [] {\includegraphics[width=7.5cm,clip]{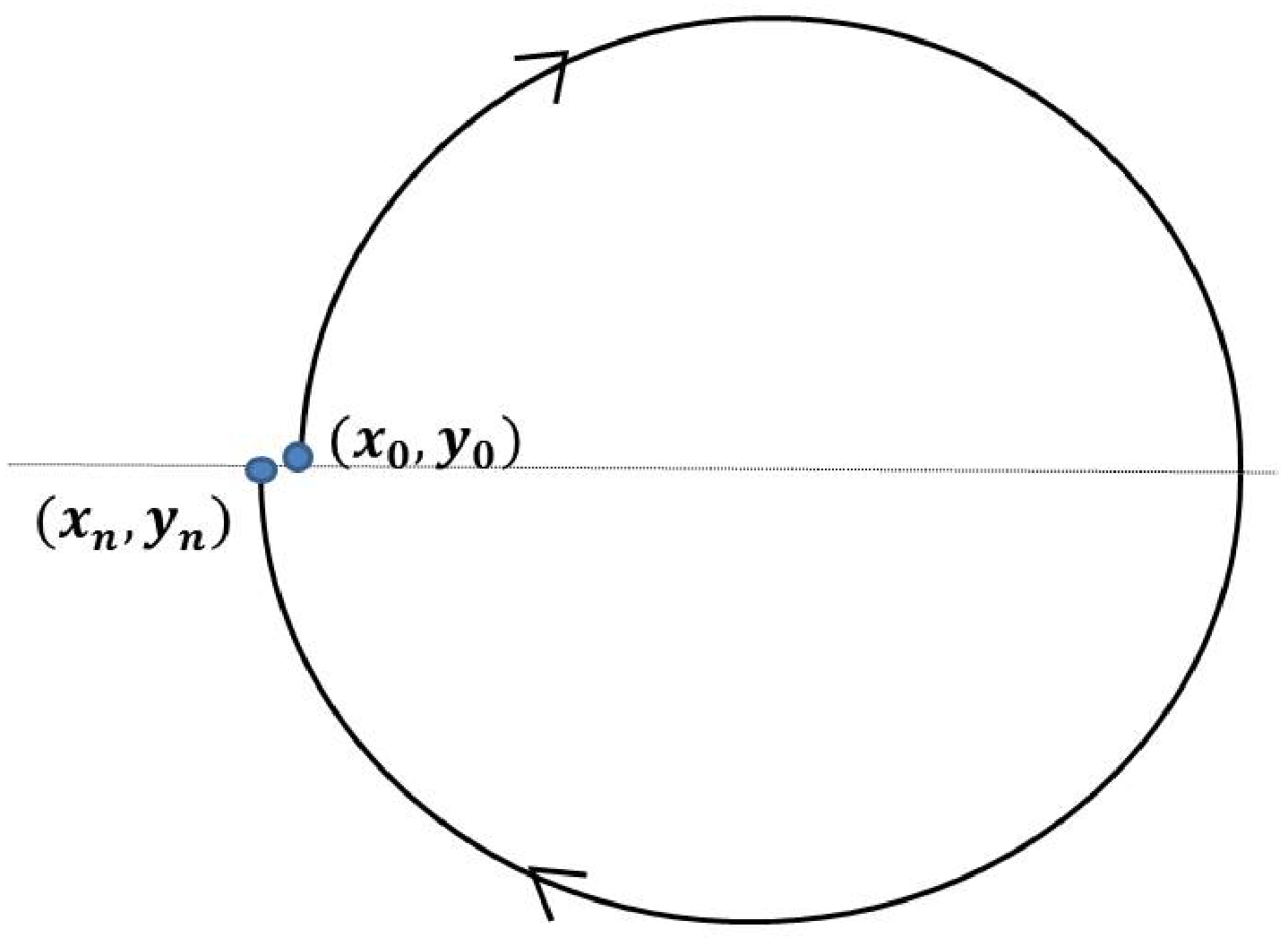} }
\subfigure [] {\includegraphics[width=7.5cm,clip]{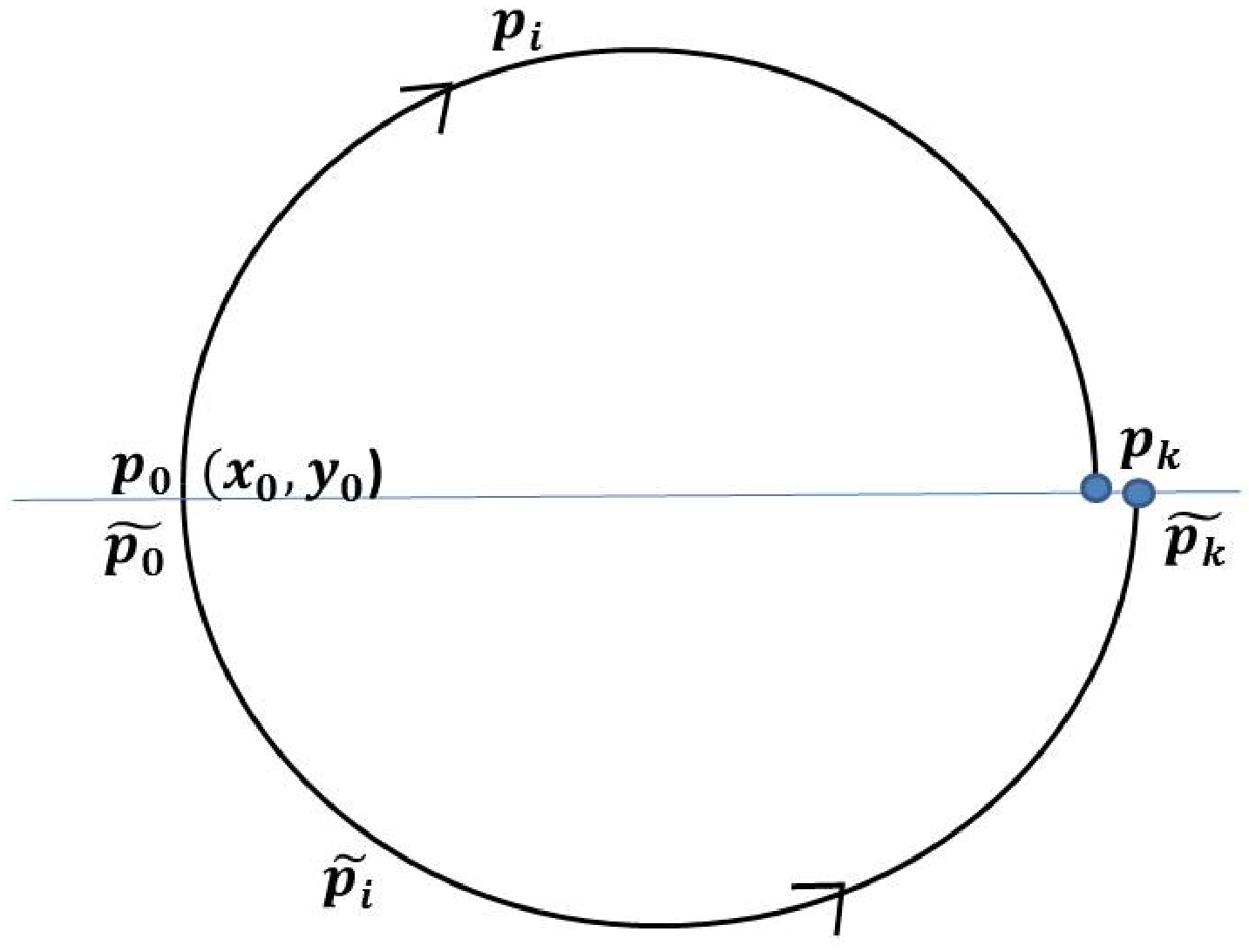} }
}
\caption{
 (a). Closed field line is solved by ``Method One".
 (b). Closed field line is solved by ``Method Two".}\label{fig:closefiled}
\end{figure}

\bigskip

\textbf{Determine the closed field lines:}
In the case of closed field lines, due to numerical errors, the approximated field line is not necessarily closed, the calculation stops when $(x_n, y_n)\in O((x_0, y_0),\delta)$ with $\delta=\min\{h_x,h_y\}$. However, when $h_x$, $h_y$ become small, there may exist no such $(x_n,y_n)$, see Figure \ref{fig:closefiled} $(a)$. In this case, the start and end points can not be considered as identical and the connection conditions in \eqref{eq:connection} can be violated. We call this approximation as ``Method One".
To avoid this problem, we can determine a closed field line by solving the following two systems simultaneously with the same time step
 \beq\label{eq:ODE2}
\left\{\begin{array}{ll}
\displaystyle \dot x(t)=\cos\theta(x(t),y(t)),\quad x(0)=x_0
\\
\displaystyle \dot y(t)=\sin\theta(x(t),y(t)),\quad y(0)=y_0
\end{array}
\right.
\quad 
\text{
and
}
\qquad
\left\{\begin{array}{ll}
\displaystyle \dot x(t)=-\cos\theta(x(t),y(t)),\quad x(0)=x_0
\\
\displaystyle \dot y(t)=-\sin\theta(x(t),y(t)),\quad y(0)=y_0
\end{array}
.\right.
\eeq
The first system in \eqref{eq:ODE2} yields a sequence of points on $\ell$, we denote them by $p_k$ ($k=0,1,\cdots$). The second system yields another sequence $\tilde p_k$ ($k=0,1,\cdots$). Then the calculations stop when $p_k\in O(\tilde{p}_k,\delta)$ with $\delta=\min\{h_x,h_y\}$ or when $\{|p_{k-1}-\tilde p_{k-1}|,|p_{k+1}-\tilde p_{k+1}|\}>|p_k-\tilde p_k|$, see Figure \ref{fig:closefiled} $(b)$.
 So the points $\{p_0,p_1,\cdots,p_k,\tilde{p}_k,\tilde p_{k-1}, \cdots,\tilde p_0\}$ form a discretization of the corresponding field line ( called it ``Method Two"). Since we will replace the local discretization at the point $p_0$ ($\tilde{p}_0$) by the field line integration, the start and end points are identical for method two, so that the interface conditions in \eqref{eq:connection} remain valid. We will see from the numerical examples that the results are different for these two different ways of determining the field lines.

\begin{figure}[htb]
\centering
{
\includegraphics[width=7.5cm] {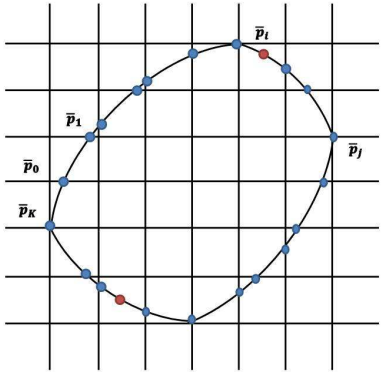}}
\caption{
{ Notations of the quadrature points along the closed  field line $\ell$. Here those points on cell edges are denoted by $\bar{p}_k$ ($k=2,3,4,\cdots,K$) and the two star points inside the cell are the discontinuities of $\nabla \cdot \mathbf{b}$.
}}\label{fig:9point}
\end{figure}

\bigskip
\textbf{Determine the quadrature points:}
The numerical quadrature points are given by interpolation.
To simplify some of the notations as well as to make the discussions clearer, we use the field line and notations in Figure \ref{fig:9point} as an example.

In order to evaluate the integral along the field line $\ell$, the quadrature points are determined by the intersection points of
 $\ell$ with those cell edges that are parallel to the $x$ axis or $y$ axis. We denote the intersection points by
 $( x_{i}, \bar y_{i})$ or $(\bar x_{j}, y_{j})$. They are numerically determined by linear interpolation of the discrete points $p_k$. It is important to note that $\nabla \cdot\mathbf{b}$ can have jumps. To get a good approximation of $E$ in \eqref{eq:E} as well as the integration in \eqref{eq:ALLt}, the discontinuities of $\nabla \cdot\mathbf{b}$  on the field line should be included in the quadrature set. We denote by $S_{\ell}=\{\bar{p}_{k}, k=0, 1, 2\cdots K_{\ell}\}$ the set of all quadrature points on the field line $\ell$.

\subsection{ 9-Point FDM }
For all those grid points inside the computational domain $\Omega\setminus\Gamma$, we use the classical $9$-Point finite difference method (FDM) to discretize the first equation in \eqref{eq:ALLt}. First of all, we use $v_{i\pm 1/2,j}\approx v(x_{i\pm 1/2},y_j)$, $v_{i,j\pm 1/2}\approx v(x_{i},y_{j\pm 1/2})$ to approximate a general function $v$. The standard centered difference approximations
read
 \beq\label{eq:9point}\begin{array}{ll}
\displaystyle \frac{\partial u^{\epsilon}}{\partial x}\Big|_{i\pm 1/2,j}\approx\frac{u_{i\pm1,j}^{\epsilon}-u_{i,j}^{\epsilon}}{\pm h_x},\quad \frac{\partial u^{\epsilon}}{\partial y}\Big|_{i,j\pm 1/2}\approx\frac{u_{i,j\pm 1}^{\epsilon}-u_{i,j}^{\epsilon}}{\pm h_y},\\
\displaystyle \frac{\partial u^{\epsilon}}{\partial y}\Big|_{i\pm 1/2,j}\approx\frac{u_{i\pm 1/2,j+1}^{\epsilon}-u_{i\pm 1/2,j-1}^{\epsilon}}{2h_y}\approx\frac{u_{i\pm 1,j+1}^{\epsilon}+u_{i,j+1}^{\epsilon}-u_{i,j-1}^{\epsilon}
-u_{i\pm 1,j-1}^{\epsilon}}{4h_y},\\
\displaystyle \frac{\partial u^{\epsilon}}{\partial x}\Big|_{i,j\pm 1/2}\approx\frac{u_{i+1,j\pm 1/2}^{\epsilon}-u_{i-1,j\pm 1/2}^{\epsilon}}{2h_x}\approx\frac{u_{i+1,j\pm 1}^{\epsilon}+u_{i+1,j}^{\epsilon}-u_{i-1,j}^{\epsilon}-u_{i-1,j\pm 1}^{\epsilon}}{4h_x}.
\end{array}
\eeq
 Let $Q=(Q^{(1)},Q^{(2)})
 =A\nabla u^{\epsilon}$. We approximate the diffusion operator at $(x_i,y_j)$ by
 \beq\label{eq:9point2}
\nabla \cdot (A\nabla u^{\epsilon})|_{i,j}=\f{\p Q^{(1)}}{\p x}\Big|_{i,j}+
\f{\p Q^{(2)}}{\p y}\Big|_{i,j}\approx\frac{Q_{i+1/2,j}^{(1)}-Q_{i-1/2,j}^{(1)}}{h_x}+\frac{Q_{i,j+1/2}^{(2)}-Q_{i,j-1/2}^{(2)}}{h_y}£¬
\eeq
with
 $$
Q_{i\pm 1/2,j} \approx A_{i\pm 1/2,j}\cdot\left(\frac{\partial u^{\epsilon}}{\partial x}\Big|_{i\pm 1/2,j},\frac{\partial u^{\epsilon}}{\partial y}\Big|_{i\pm 1/2,j}\right)^T,\qquad
Q_{i,j\pm 1/2} \approx A_{i,j\pm 1/2}\cdot\left(\frac{\partial u^{\epsilon}}{\partial x}\Big|_{i,j\pm 1/2},\frac{\partial u^{\epsilon}}{\partial y}\Big|_{i,j\pm 1/2}\right)^T.
 $$
The boundary values are given by the Dirichlet boundary conditions.

 The major difference and difficulty is about discretization of the integration in the third equation in \eqref{eq:ALLt}.  The calculation is separated into two steps:
 \begin{itemize}
 \item
Calculate $E(s)$. The the quadrature points ${\bar p_{i}}( i=0,1,\cdots K_{\ell})$ on each field line $\ell$ are determined as in section 3.1.
 We use the composite trapezoidal rule to approximate the integration in $E(s)=e^{\int_0^ {{s}} \nabla \cdot \mathbf{b}  \,\,ds'}+e^{-\int_{L_{\ell}}^ {{s}} \nabla \cdot \mathbf{b} \,\,ds'}$. Since the discontinuities in $\nabla\cdot\mathbf{b}$ have been taken into account in the quadrature set, we calculate the value of $E(s)$ at the point $p_k$ by
\beq\label{eq:discE}\exp\Big(\sum_{i=0}^{k-1}\f{\omega_i}{2} \big(\nabla\cdot \mathbf{b}|_{\bar p_{i}}+\nabla \cdot\mathbf{b}|_{\bar p_{i+1}}\big)\Big)+\exp\Big(-\sum_{i=k+1}^{K_\ell}\f{\omega_{i-1}}{2} \big(\nabla \cdot\mathbf{b}|_{\bar p_{i-1}}+\nabla \cdot\mathbf{b}|_{\bar p_{i}}\big)\Big).
    \eeq
Here $\omega_i=|\bar p_{i+1}-\bar p_i|$ is the Eulerian distance between the two quadrature points.

  \item
The diffusion operator $- \nabla
   \cdot \big( \alpha \mathbf{b}_{\perp} (\mathbf{b}_{\perp}\cdot \nabla u^{\epsilon})\big)$ on each grid point is approximated by the centered finite difference method same as \eqref{eq:9point2}. Let $S_1$, $S_2$ be respectively the sets of quadrature points on the edge parallel to the $x$-axis and $y$-axis. The value of $\Theta=f+\nabla
   \cdot \big( \alpha \mathbf{b}_{\perp} (\mathbf{b}_{\perp}\cdot \nabla u^{\epsilon})\big)$ at the quadrature points
   $p_k$ can be given by a linear interpolation such that:
   \beq
  \Theta\mid_{\bar p_k}\approx\left\{\begin{array}{ll} \Theta\mid_{(x_{i},y_{j})}\frac{x_i+h_x-\bar x_i}{h_x}+
  \Theta\mid_{(x_{i}+h_x,y_{j})}\frac{\bar x_i-x_i}{h_x},\quad  \mbox{for }\bar p_k=(\bar x_i, y_j)\in S_{1}, \bar x_i\in(x_i, x_i+h_x),
\\ \Theta\mid_{(x_{i},y_{j})}\frac{y_j+h_y-\bar y_j}{h_y}+
  \Theta\mid_{(x_{i},y_{j}+h_y)}\frac{\bar y_j-y_j}{h_y},\quad  \mbox{for }\bar p_k=(x_i, \bar y_j)\in S_{2}, \bar y_j\in(y_j, y_j+h_y).
    \end{array}
    \right.
  \eeq
The integration in the third equation of \eqref{eq:ALLt} is approximated by
\beq\label{eq:discL}\begin{aligned} &\sum_{i=0}^{K_{\ell}-1}\f{\omega_i}{2} \Big((E\Theta)\mid_{\bar p_i}+(E\Theta)\mid_{\bar p_{i+1}}\Big)
    =0,
\end{aligned}\eeq
with
 $
\omega_i=
|p_{i+1}-p_i|$.

For $\bar p_k=(x^\ast,y^\ast)$ inside the cell $[x_i,x_{i+1}]\times[y_j,y_{j+1}]$, $\Theta\mid_{\bar p_k}$ is given by the following bilinear interpolation
$\Theta\mid_{\bar p_k}=\varpi_1 \Theta\mid_{(x_i, y_j)}$ $+\varpi_2 \Theta\mid_{(x_{i+1}, y_j)}$ $+\varpi_3 \Theta\mid_{(x_i, y_{j+1})}$ $+\varpi_4 \Theta\mid_{(x_{i+1}, y_{j+1})}$ with $$\begin{array}{ll}\varpi_1=\f{(y_{j+1}-y^\ast)(x_{i+1}-x^\ast)}{h_xh_y},\quad &
    \varpi_2=\f{(y_{j+1}-y^\ast)(x^\ast-x_i)}{h_xh_y},\\[6pt] 
		\varpi_3=\f{(y^\ast-y_{j})(x_{i+1}-x^\ast)}{h_xh_y},\quad &\varpi_4=\f{(y^\ast-y_{j})(x^\ast-x_i)}{h_xh_y}.\end{array}$$

\end{itemize}


\section{Numerical Results}

We present several tests to demonstrate the performance of the proposed scheme. For all examples, the precise discretizations in section 3 are used.

\begin{figure}[htb]
\centering
{
\subfigure [] {\includegraphics[width=7.5cm,clip]{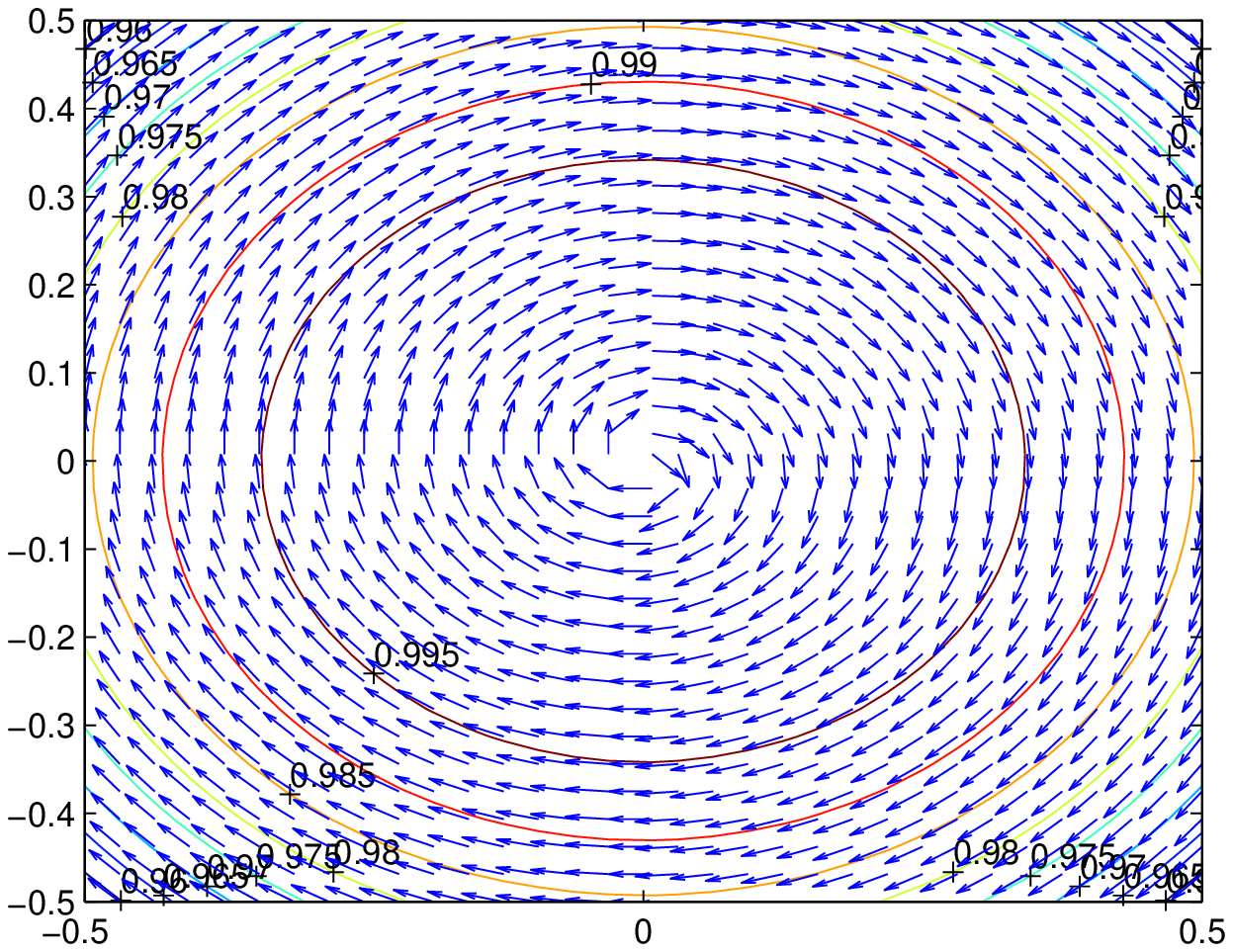} }
\subfigure [] {\includegraphics[width=7.5cm,clip]{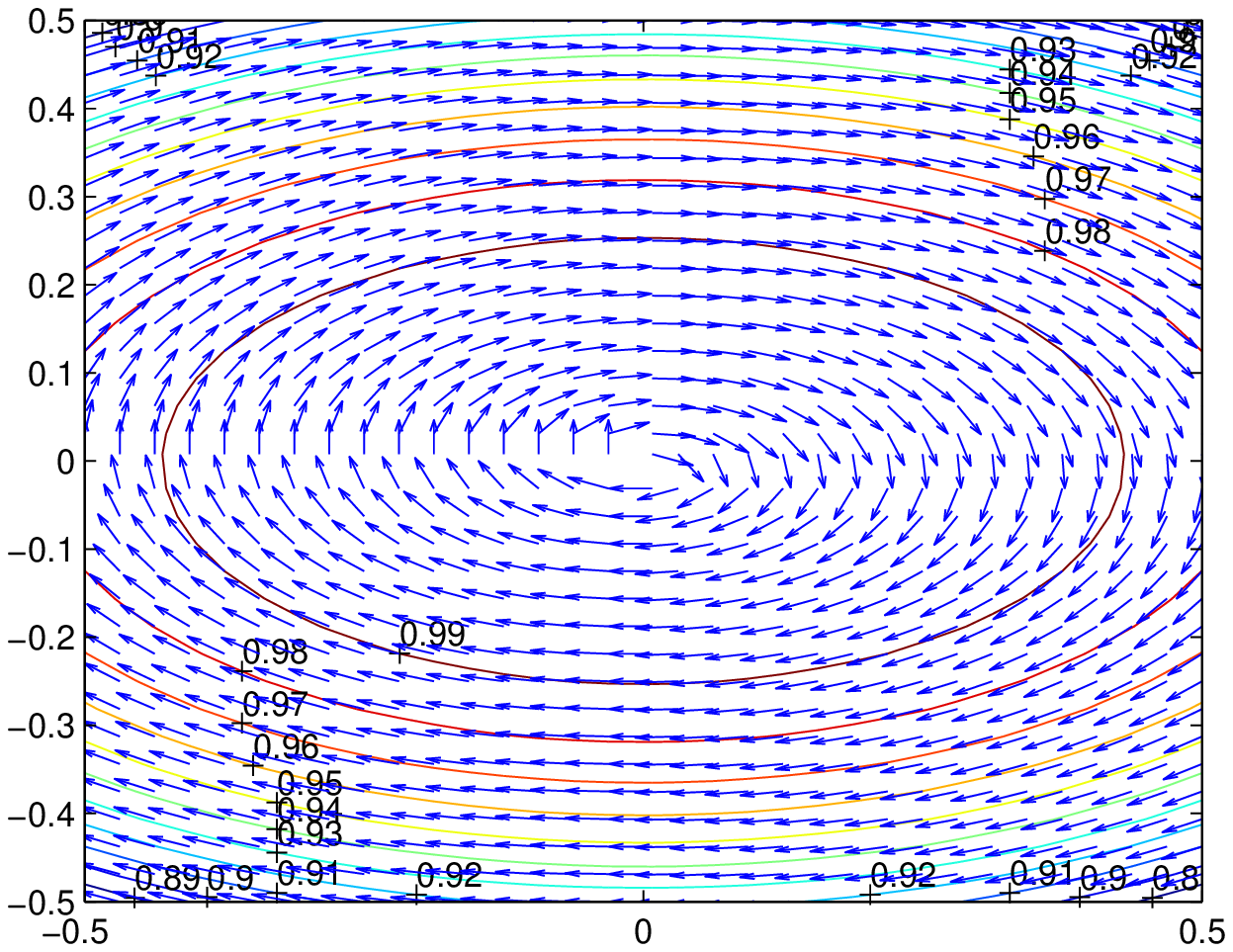} }
\subfigure[] {\includegraphics[width=7.5cm,clip]{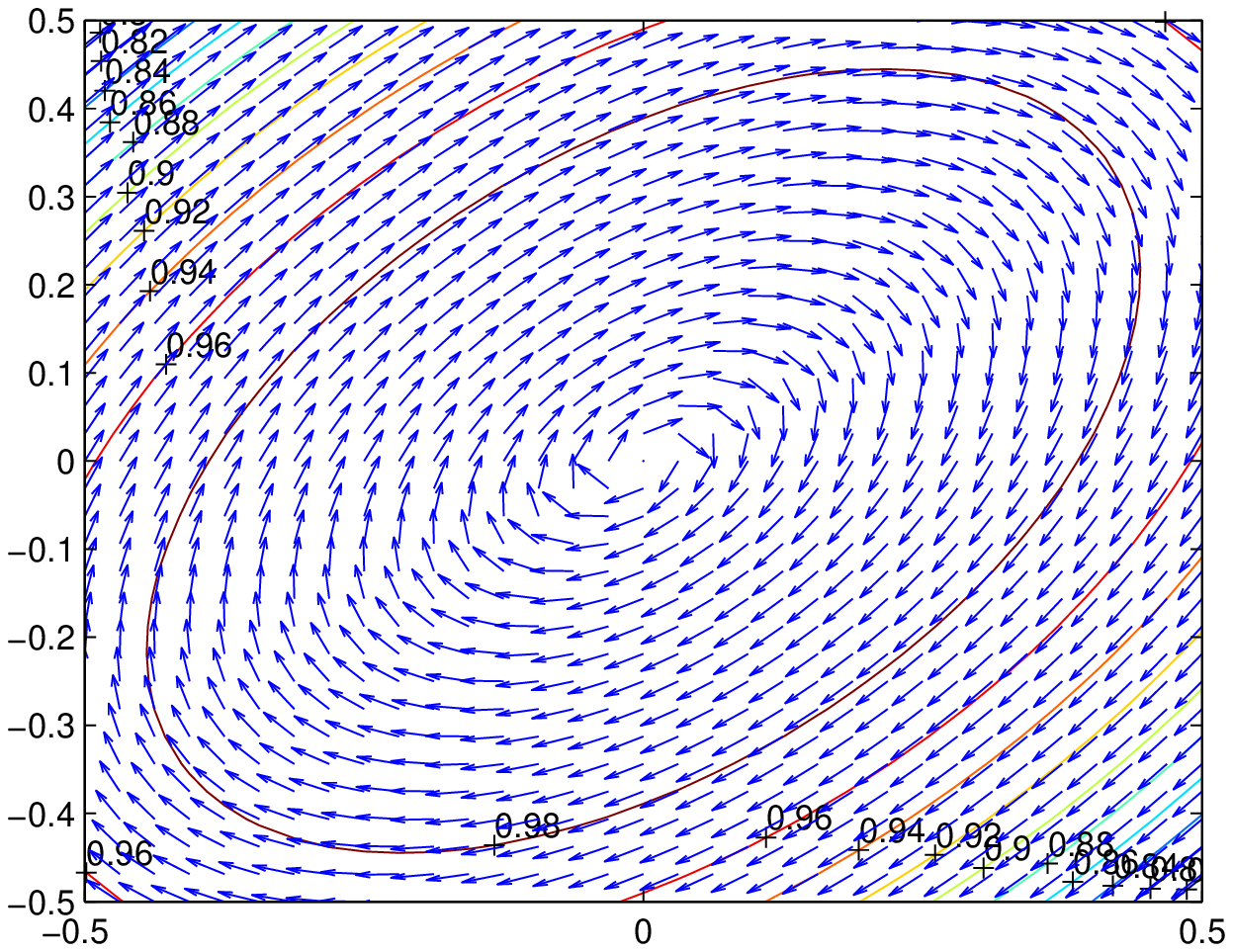}  }
\subfigure [] {\includegraphics[width=7.5cm,clip]{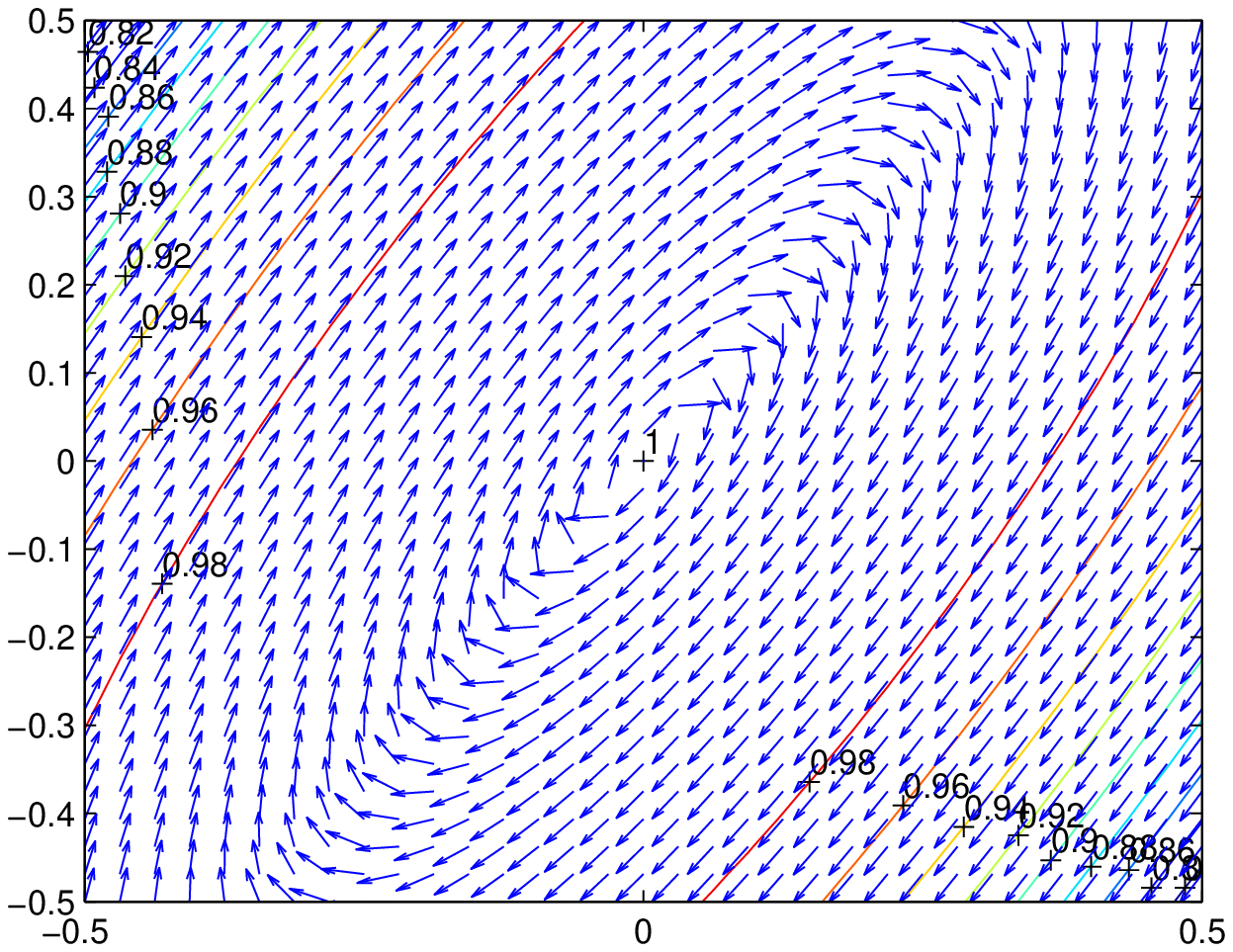} }

}
\caption{Example 1. The fields of the four test cases.
 (a). $\gamma_1=\gamma_2=0.5,~\varphi=0$.
 (b). $\gamma_1=0.5,~\gamma_2=0.85,~\varphi=0$.
 (c). $\gamma_1=0.5,~\gamma_2=0.85,~\varphi=\pi/4$.
 (d). $\gamma_1=0.25,~\gamma_2=0.85,~\varphi=\pi/3$.}\label{fig:ex1field}
\end{figure}

\bigskip
\textbf{Example 1:}  In this example, the following  exact solution for $(x,y)\in[-0.5,0.5]\times[-0.5,0.5]$ is considered:
\beq\label{eq:uex}
u(x,y)=1-[\gamma_1^2(x\cos\varphi+y\sin\varphi)^2+\gamma_2^2(x\sin\varphi-y\cos\varphi)^2]^{3/2}.
\eeq
We test four different choices of $\gamma_1$, $\gamma_2$ and $\varphi$. All cases include closed field lines (see Figure \ref{fig:ex1field}). In the first two cases, the closed field lines are ellipses that are
symmetric with respect to the $x$ and $y$ axises, while in the last two cases, the ellipses are rotated.

For the first symmetric cases, there exist schemes in the literature \cite{Bram 14,Gunter07,Gunter05} that exhibit uniform second order convergence with respect to $\epsilon$. However,
 when the closed field lines become tilted ellipses,
 the convergence order of all schemes discussed in \cite{Bram 14,Gunter07,Gunter05}  depends on the level of anisotropy. In particular, when $\epsilon$ reaches the order of $10^{-9}$, no convergence can be observed with $h_x,h_y\gg \epsilon$ in the tilted ellipse case.

The field and diffusion tensor are given by
$$
 \mathbf{b}=\frac{1}{\sqrt{u_x^2+u_y^2}}\left(\begin{array}{c}-u_y\\u_x\end{array}\right)=\left(\begin{array}{c}b_1\\b_2\end{array}\right),
\qquad A=\left(\begin{array}{cc}b_1&-b_2\\b_2&b_1\end{array}\right)\left(\begin{array}{cc}1/\epsilon &0\\0&1\end{array}\right)
\left(\begin{array}{cc}b_1&b_2\\-b_2&b_1\end{array}\right),
$$
 while the boundary conditions and the source term are determined by substituting the exact solutions into \eqref{eq:ellipT}. At the origin $(0,0)$, $b$ has no definition. We replace ${\sqrt{u_x^2+u_y^2}}$ by ${\sqrt{u_x^2+u_y^2+\delta}}$ $(\delta=O(10^{-16}))$ to approximate the solution at the origin.

\begin{figure}[htb]
\centering
{
\subfigure [] {\includegraphics[width=5cm,clip]{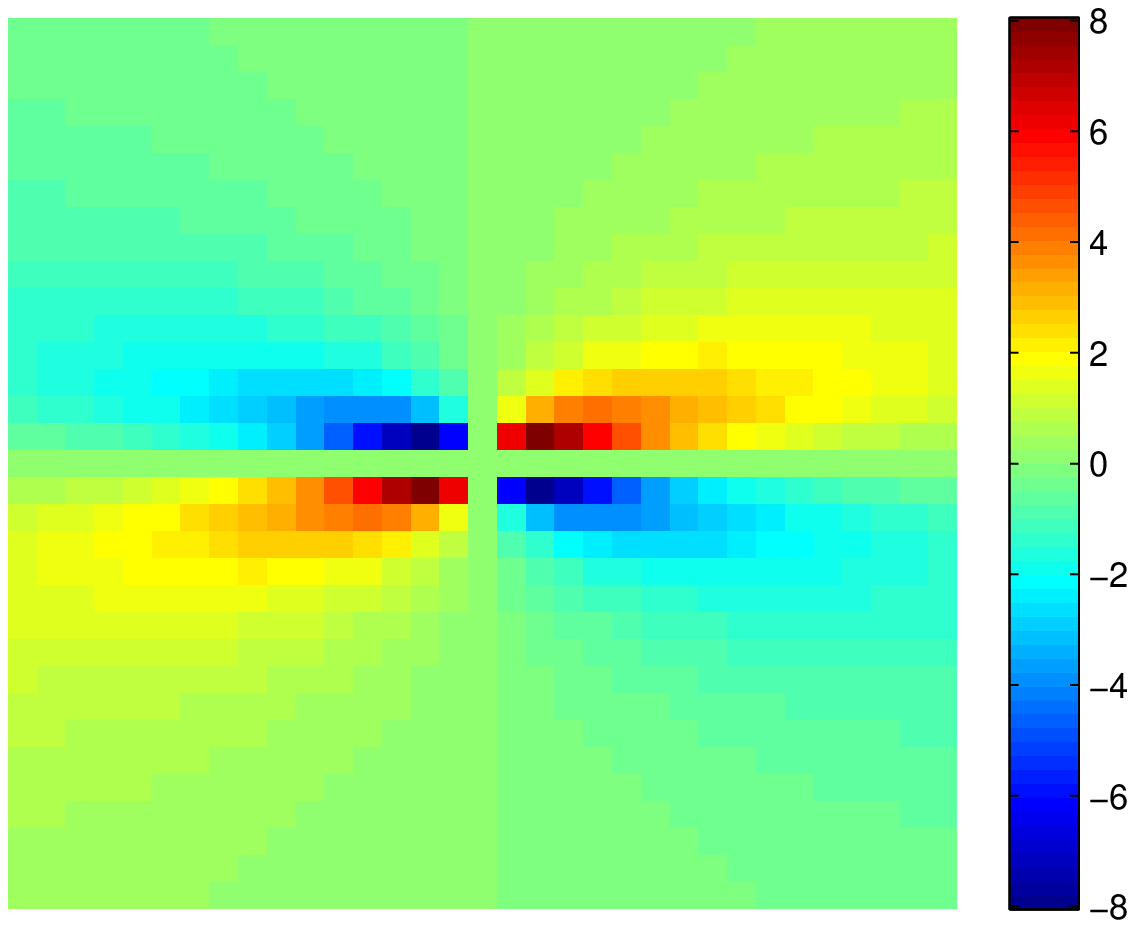} }
\subfigure [] {\includegraphics[width=5cm,clip]{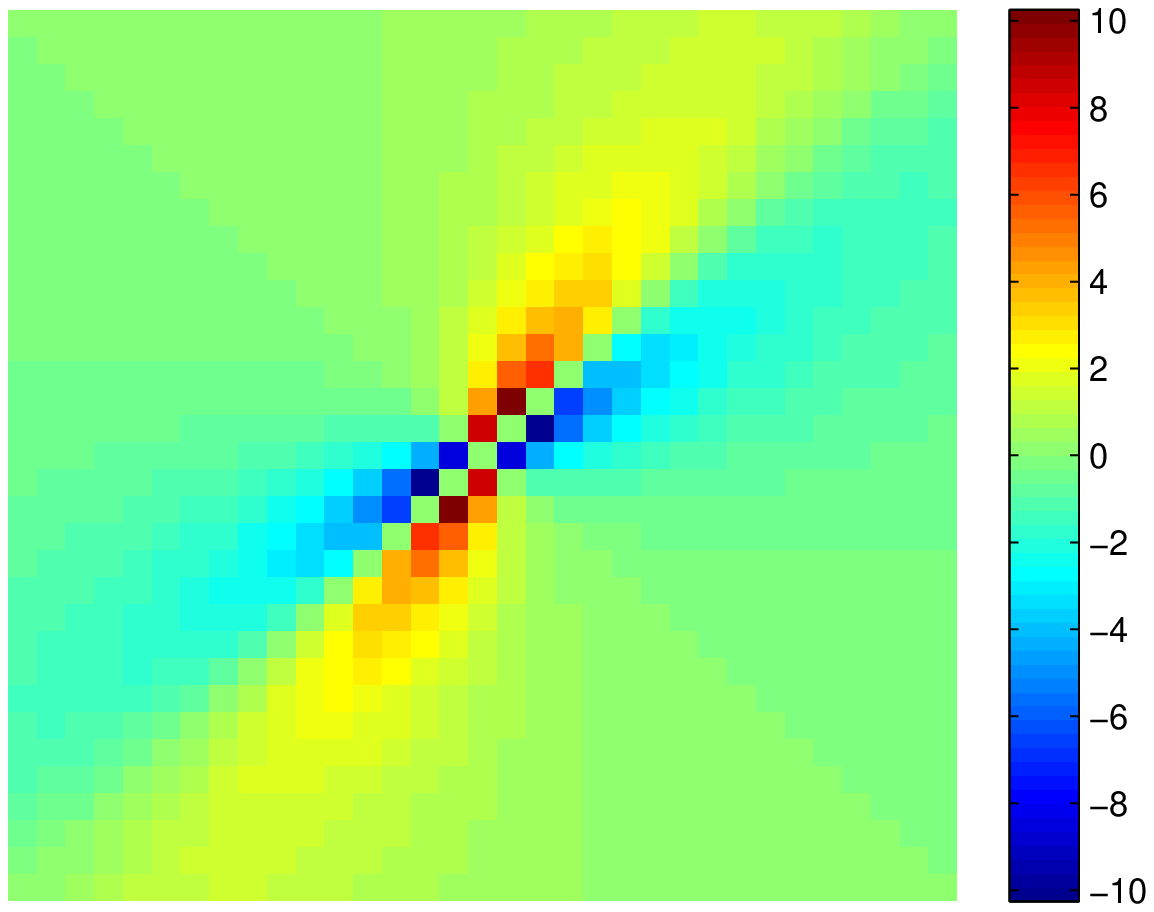} }
\subfigure[] {\includegraphics[width=5cm,clip]{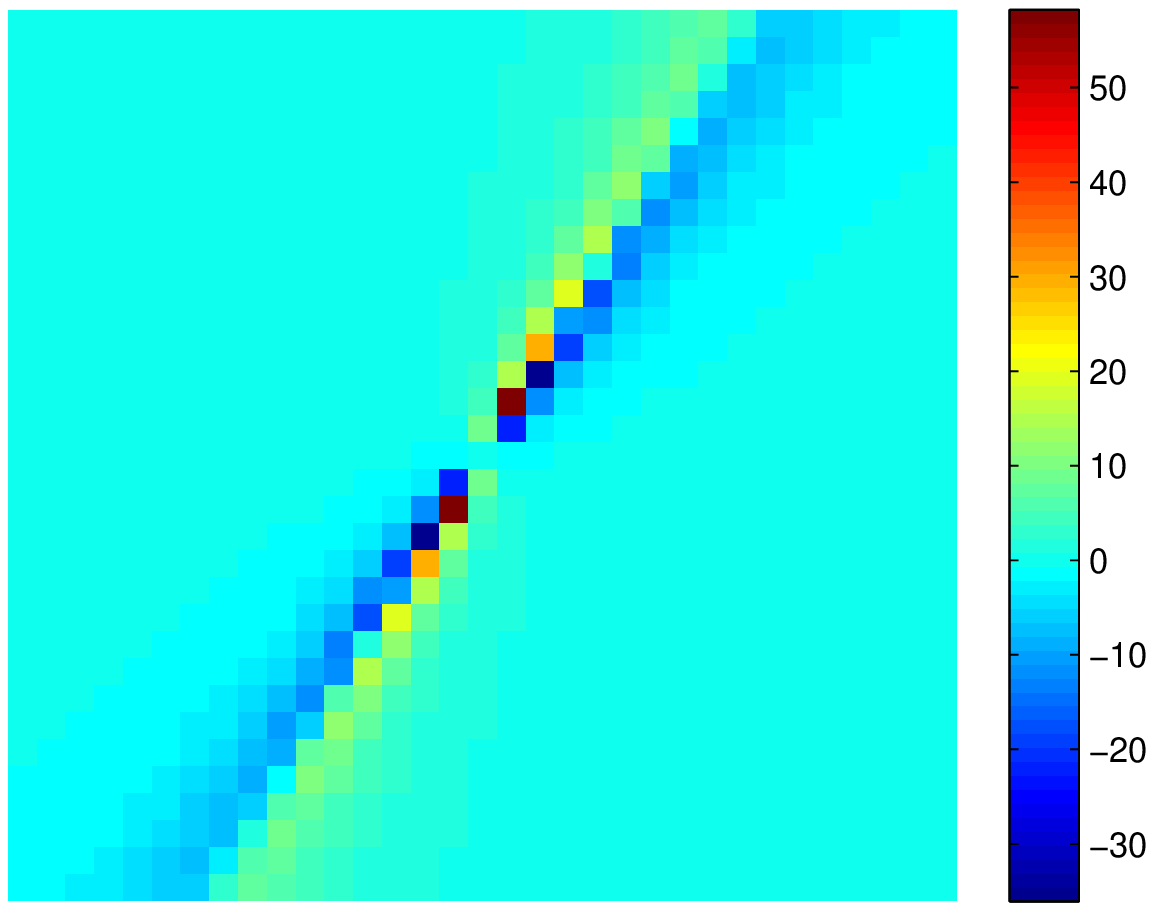}  }

}
\caption{Example 1. $\nabla\cdot\mathbf{b}$, $I\times J=32\times 32$.
 (a). $\gamma_1=0.5,~\gamma_2=0.85,~\varphi=0$.
 (b). $\gamma_1=0.5,~\gamma_2=0.85,~\varphi=\pi/4$.
 (c). $\gamma_1=0.25,~\gamma_2=0.85,~\varphi=\pi/3$.}\label{fig:ex1divb}
\end{figure}

In the first case, $\nabla \cdot\mathbf{b}=0$ on each field line. In all other cases, $\nabla \cdot\mathbf{b}\neq 0$, and the discontinuities of $\nabla \cdot\mathbf{b}$ are distributed on a ray. The angle between the ray and the x-axis is equal to $\varphi$, as shown in Figure \ref{fig:ex1divb}.

\begin{figure}[htb]
\centering
{
\subfigure [] {\includegraphics[width=6cm,clip]{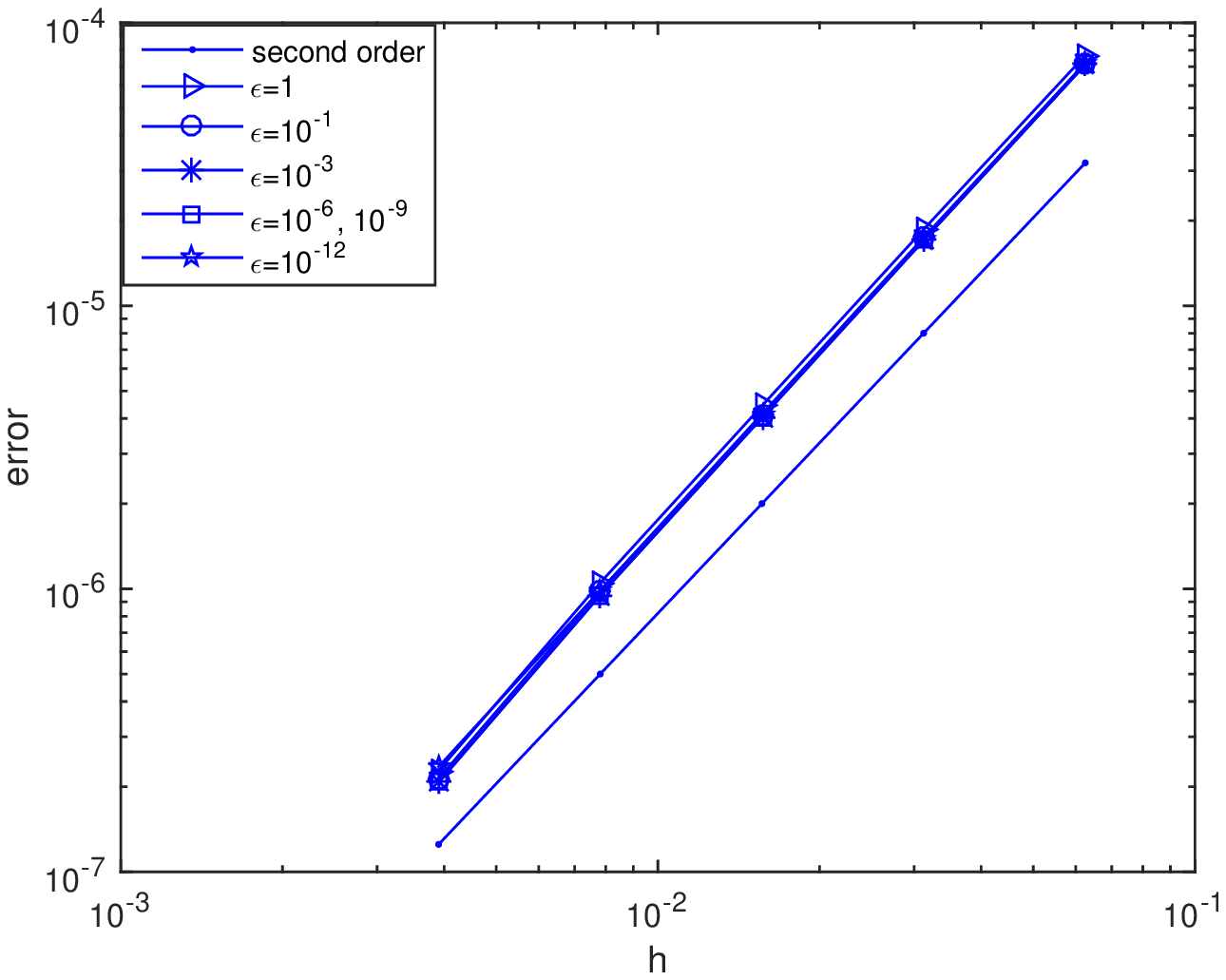} }
\subfigure [] {\includegraphics[width=6cm,clip]{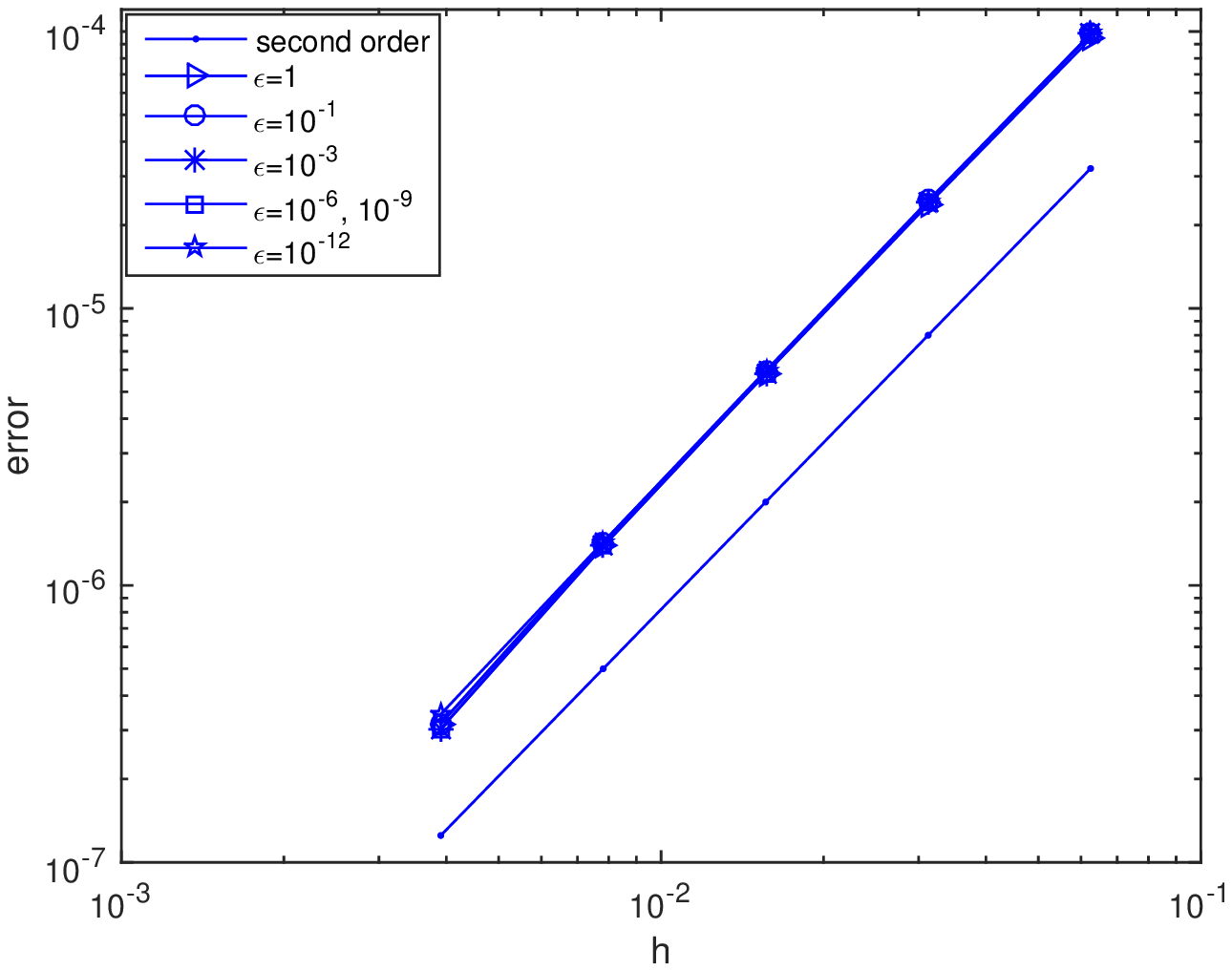} }
\subfigure[] {\includegraphics[width=6cm,clip]{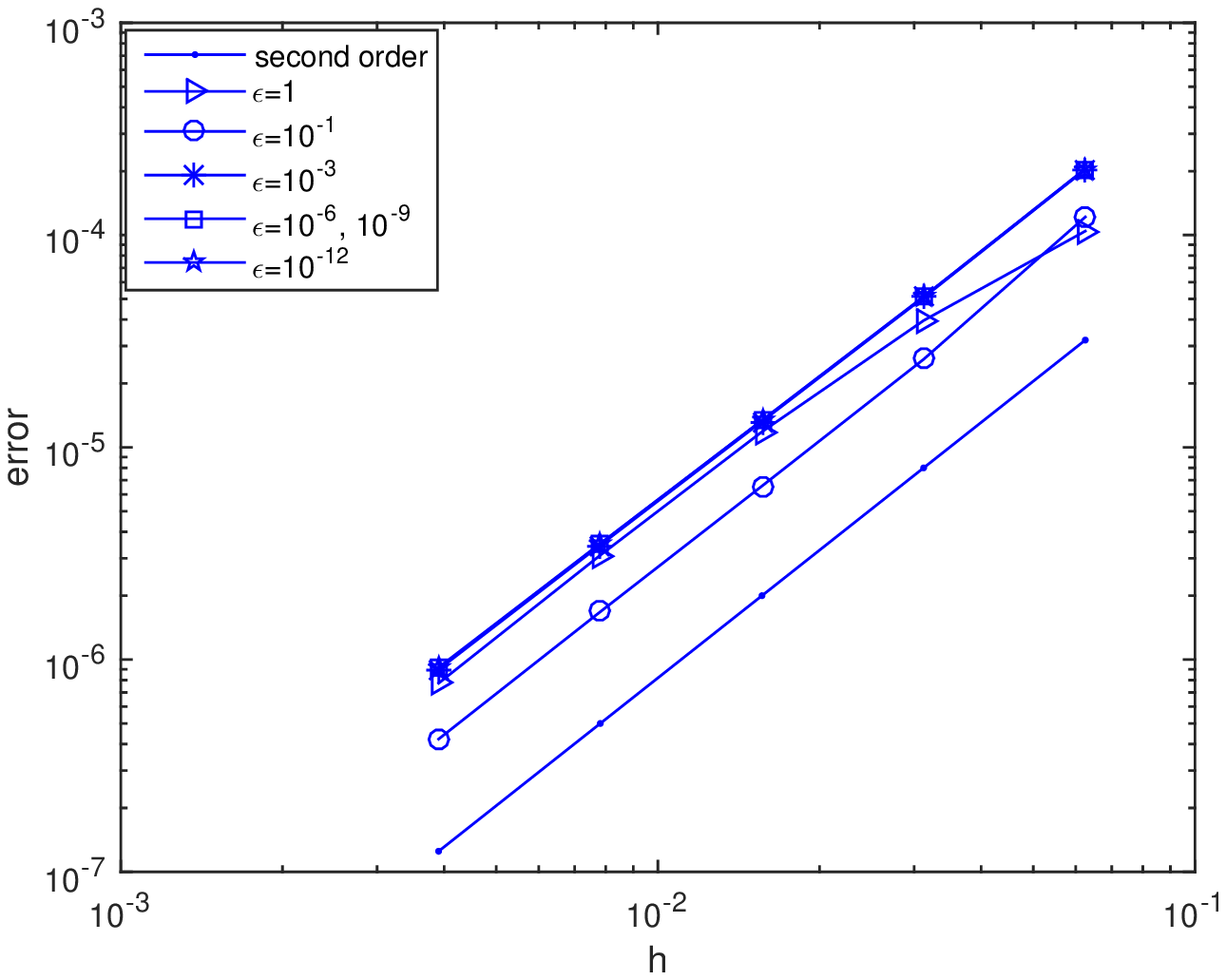}  }
\subfigure [] {\includegraphics[width=6cm,clip]{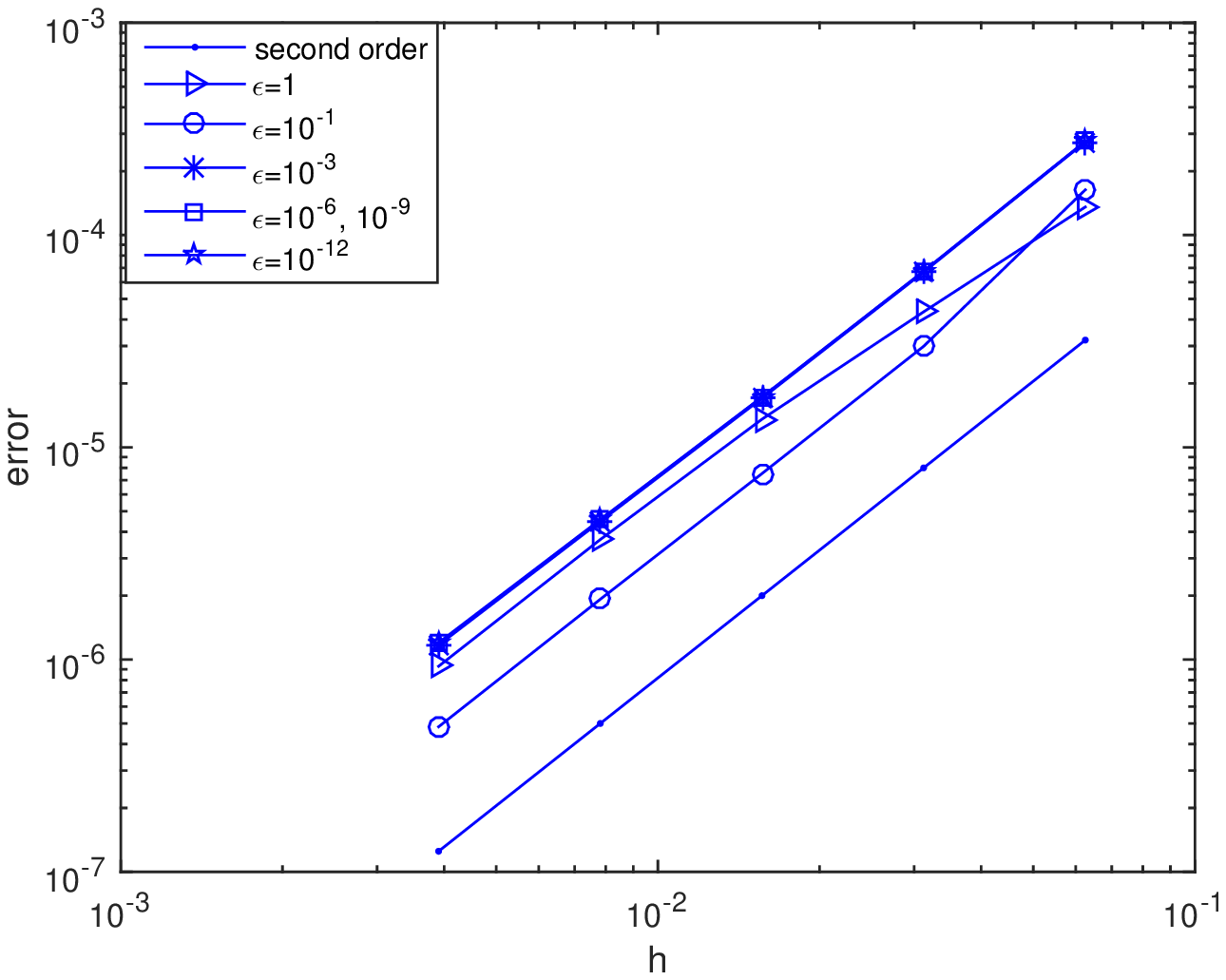} }
\subfigure [] {\includegraphics[width=6cm,clip]{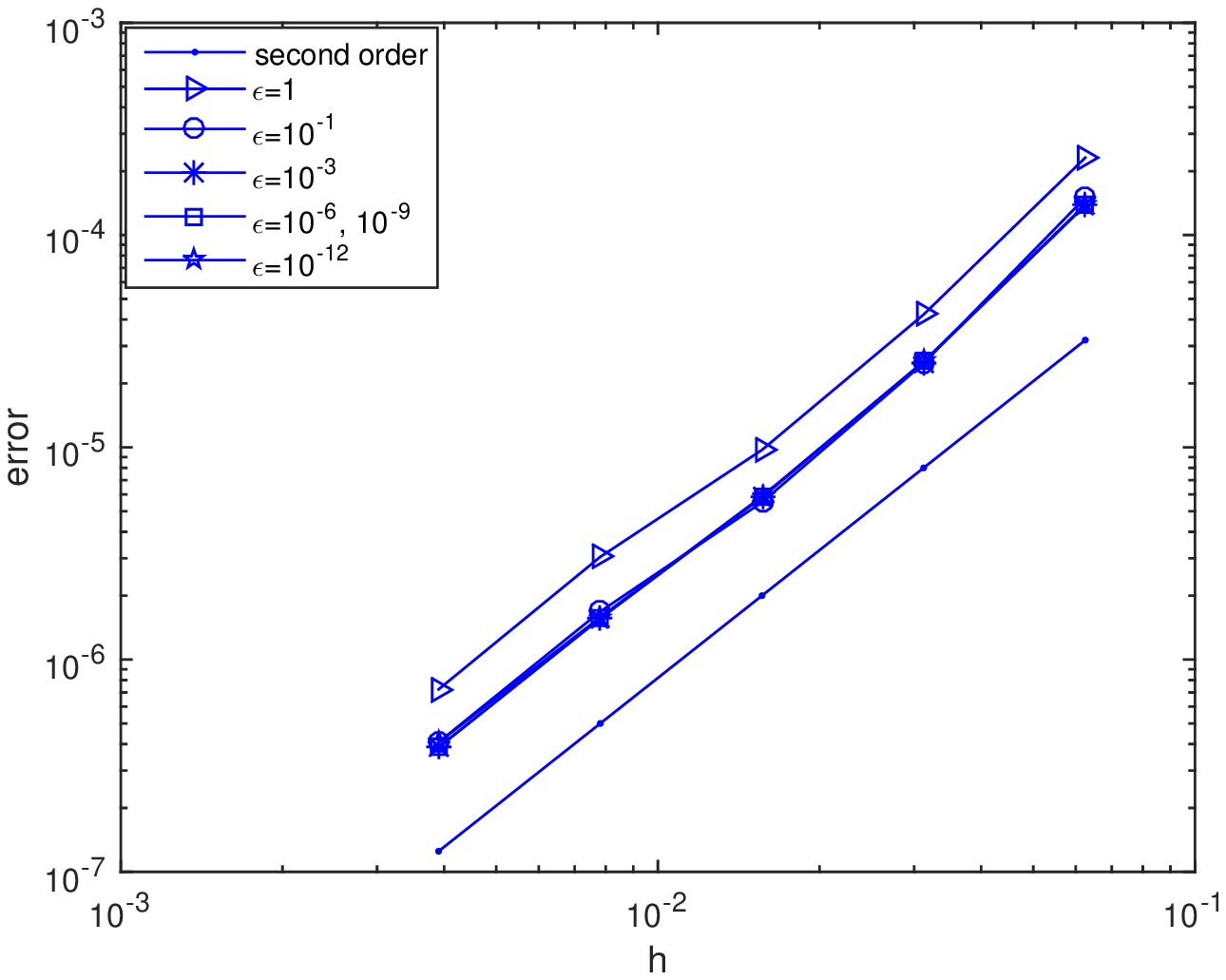} }
\subfigure [] {\includegraphics[width=6cm,clip]{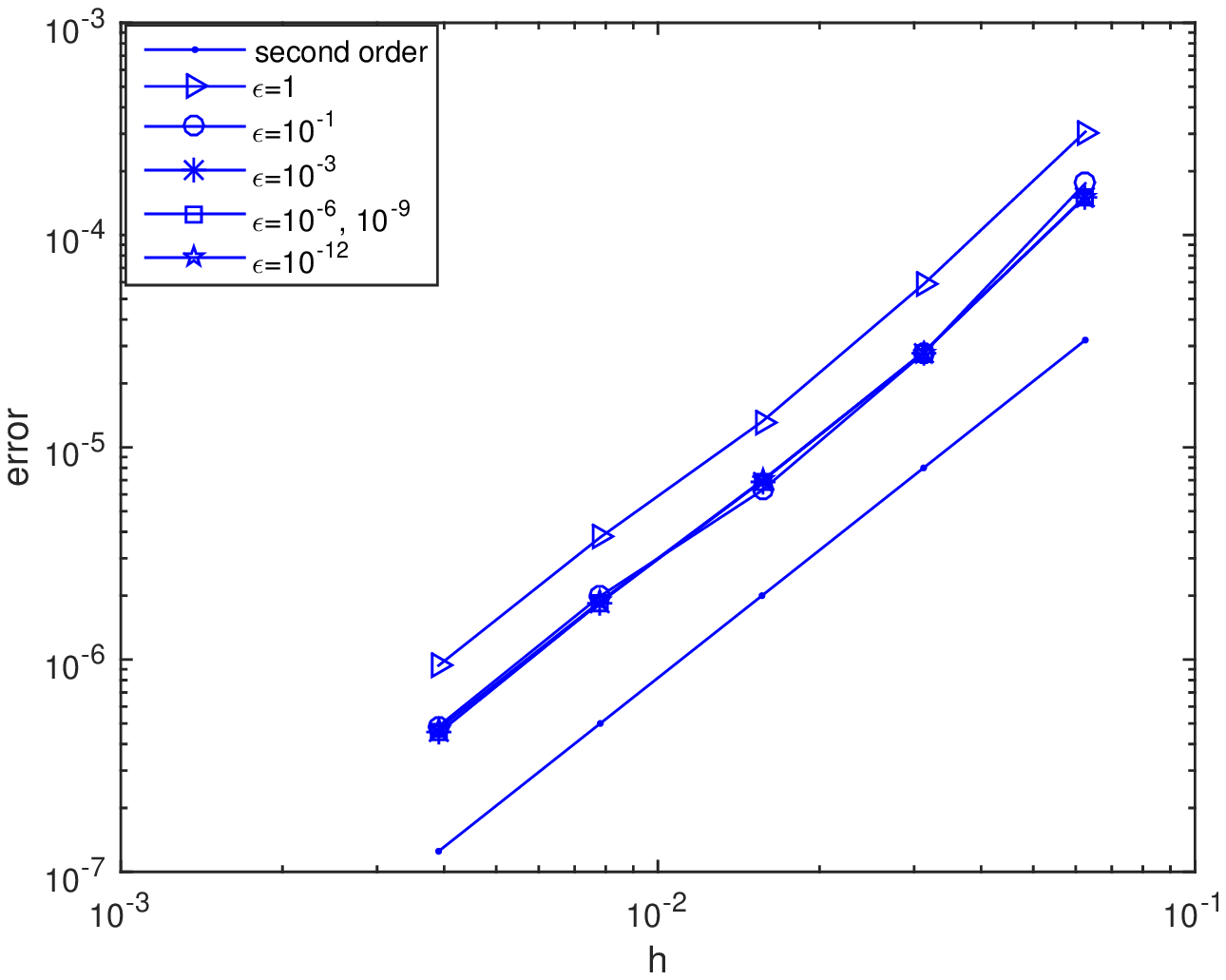} }
\subfigure[] {\includegraphics[width=6cm,clip]{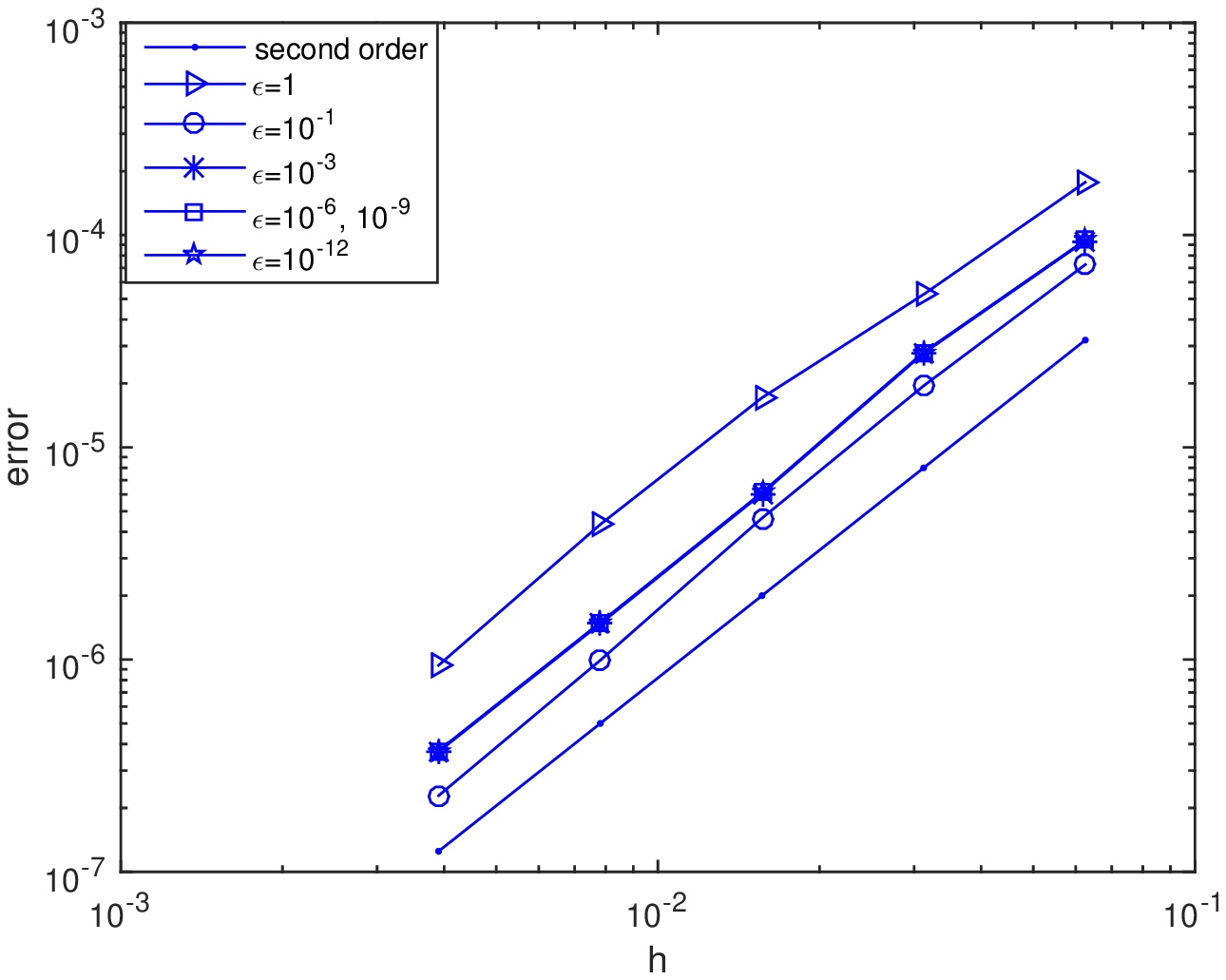}  }
\subfigure [] {\includegraphics[width=6cm,clip]{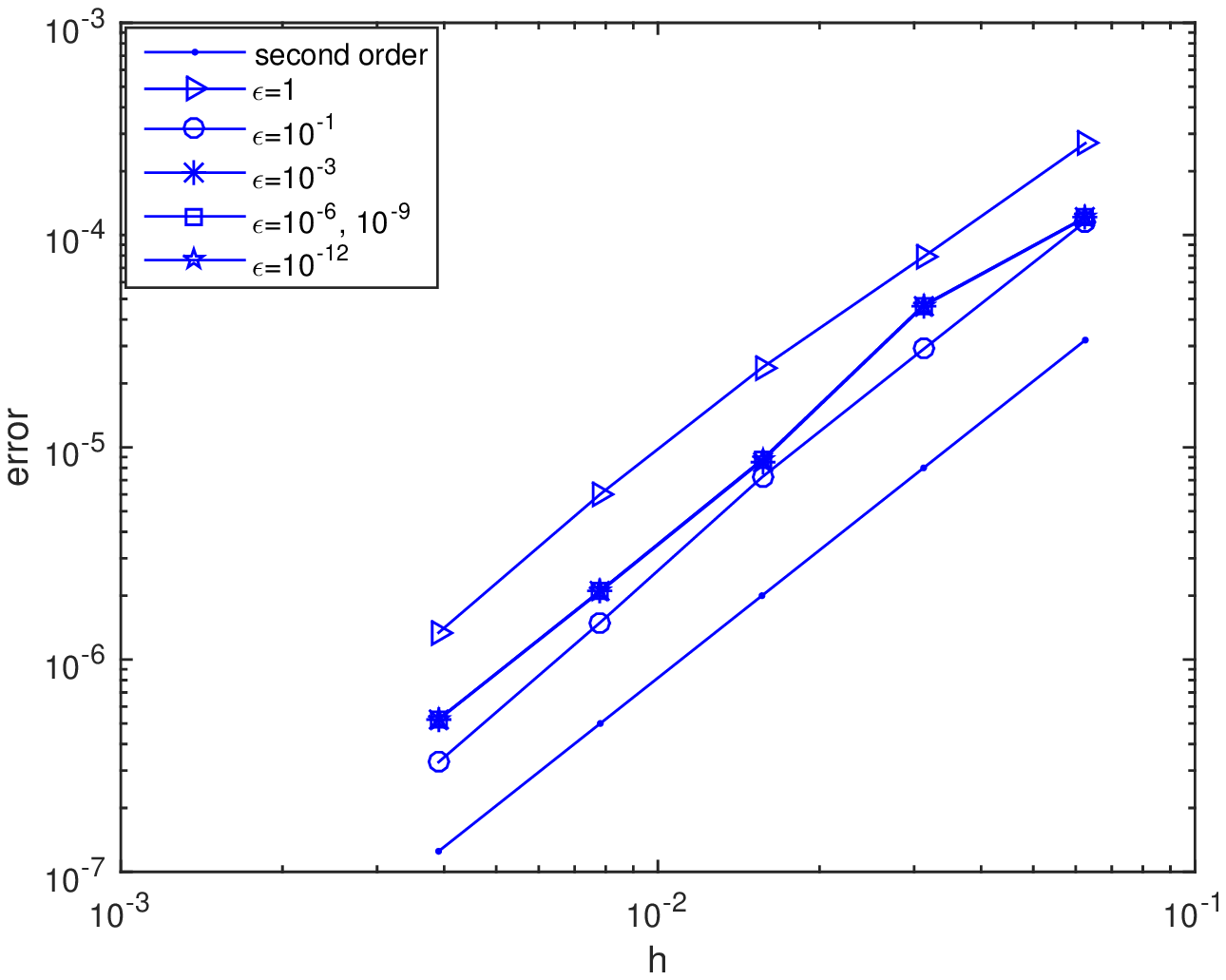} }

}
\caption{Example 1. Convergence orders for $L^{2}$ error (Left) and $L^{\infty}$ error (Right).
 (a). (b). $\gamma_1=\gamma_2=0.5,~\varphi=0$.
 (c). (d). $\gamma_1=0.5,~\gamma_2=0.85,~\varphi=0$.
 (e). (f). $\gamma_1=0.5,~\gamma_2=0.85,~\varphi=\pi/4$.
 (g). (h). $\gamma_1=0.25,~\gamma_2=0.85,~\varphi=\pi/3$.}\label{fig:ex1norm}
\end{figure}

\begin{figure}[htb]
\centering
{
\subfigure [] {\includegraphics[width=7cm,clip]{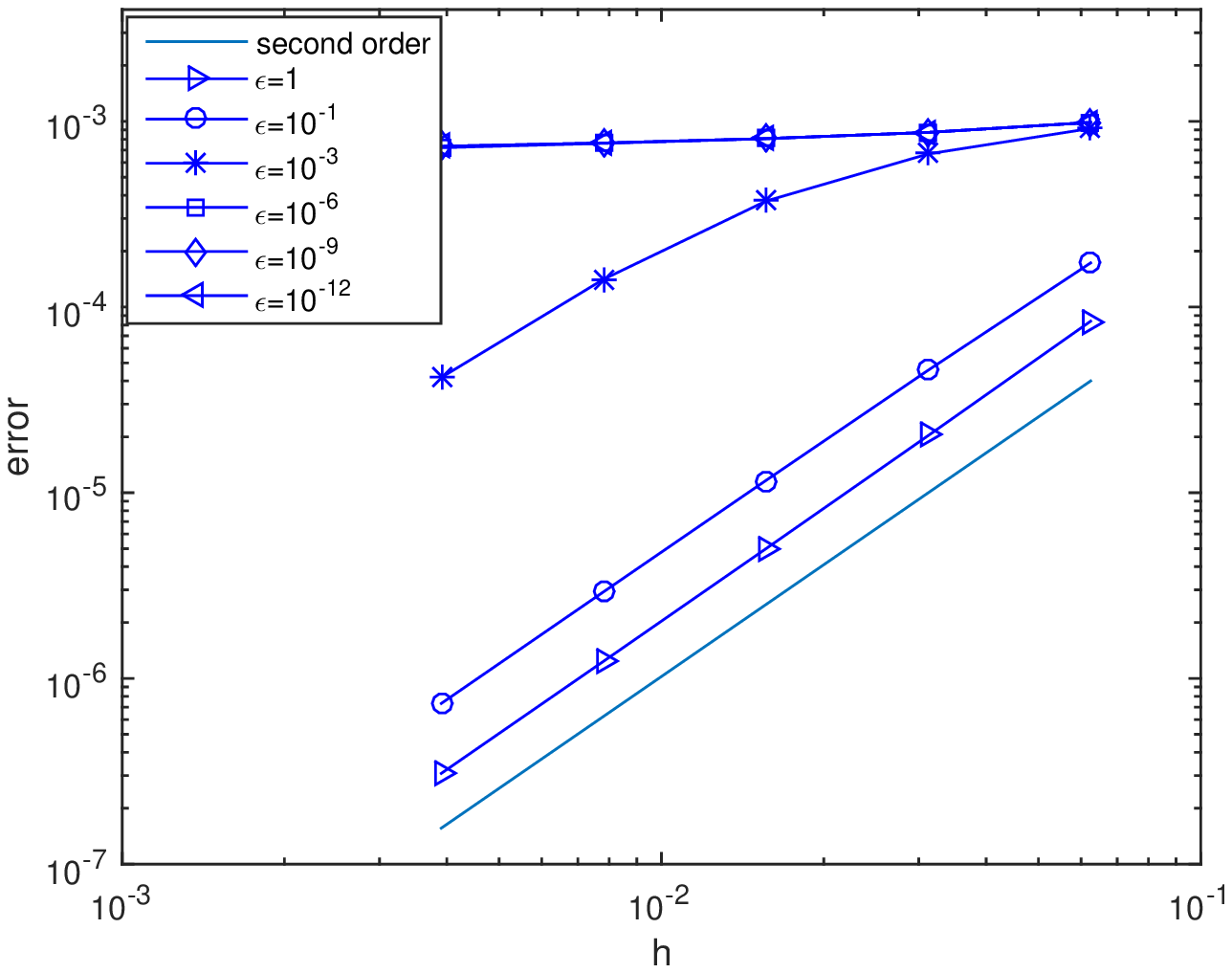} }
\subfigure[] {\includegraphics[width=7cm,clip]{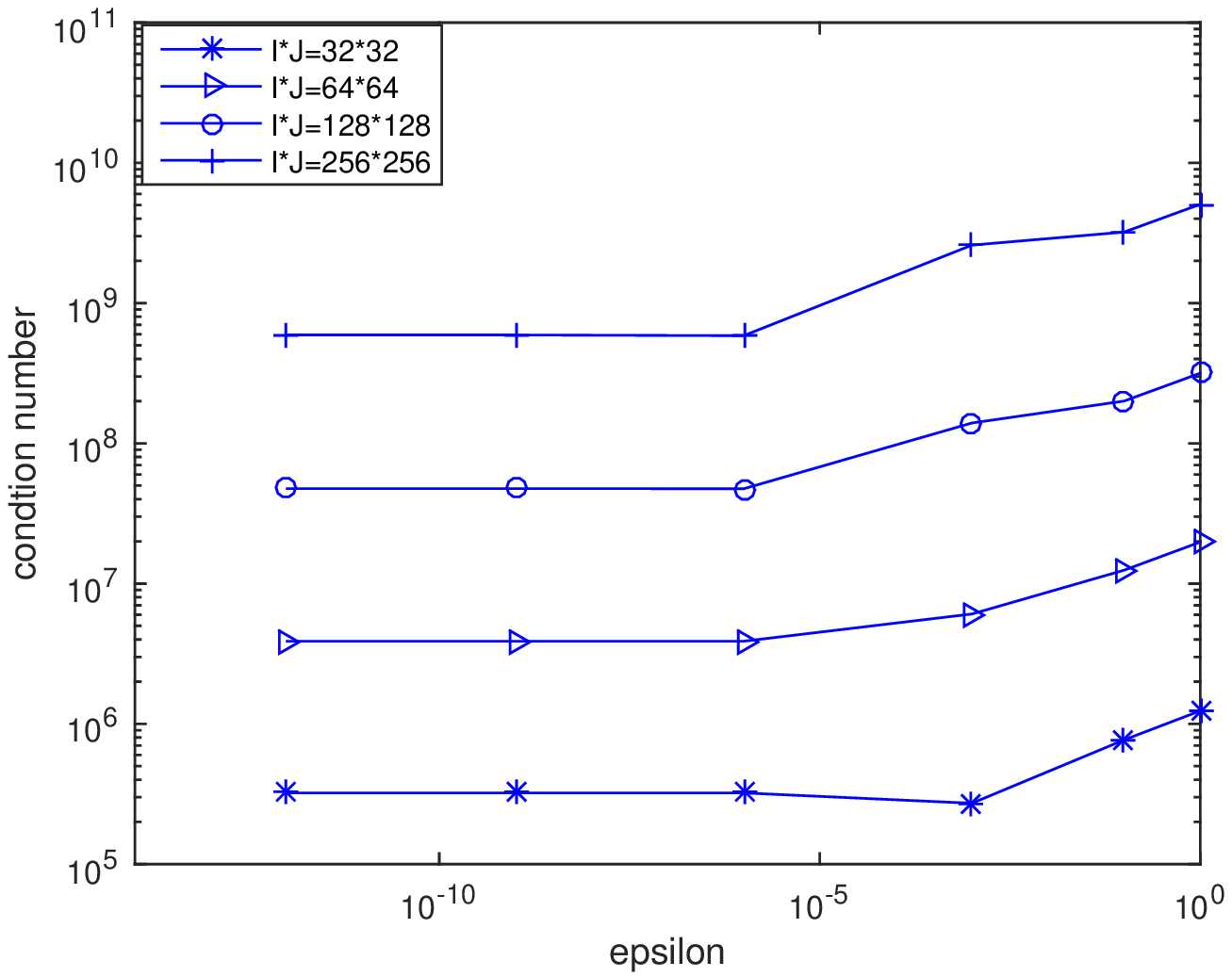}  }

}
\caption{Example 1. $\gamma_1=0.5,~\gamma_2=0.5,~\varphi=0$. Convergence orders  for $L^{\infty}$ error and Condition number estimate for the classical $9$-point FDM. (a). Convergence orders.
 (b). Condition number.}\label{fig:ex1d1}
\end{figure}

\textbf{Convergence Order }
The convergence orders of our new scheme are displayed in Figure \ref{fig:ex1norm}. Second order convergence in both $L^2$ and $L^\infty$ norm can be observed regardless of the anisotropy strength. Figure \ref{fig:ex1d1} $(a)$ shows that, when $\epsilon$ is less than the order of $10^{-3}$, the numerical errors by the classical $9$-Point FDM of the original system \eqref{eq:ellipT} does not decrease as the mesh is refined. 
When the closed field lines become tilted ellipses, the convergence order of all schemes discussed in \cite{Bram 14,Gunter07,Gunter05} depends on the level of anisotropy. The comparison illustrates the advantage of our reformulation.

\begin{figure}[htb]
\centering
{
\subfigure [] {\includegraphics[width=7cm,clip]{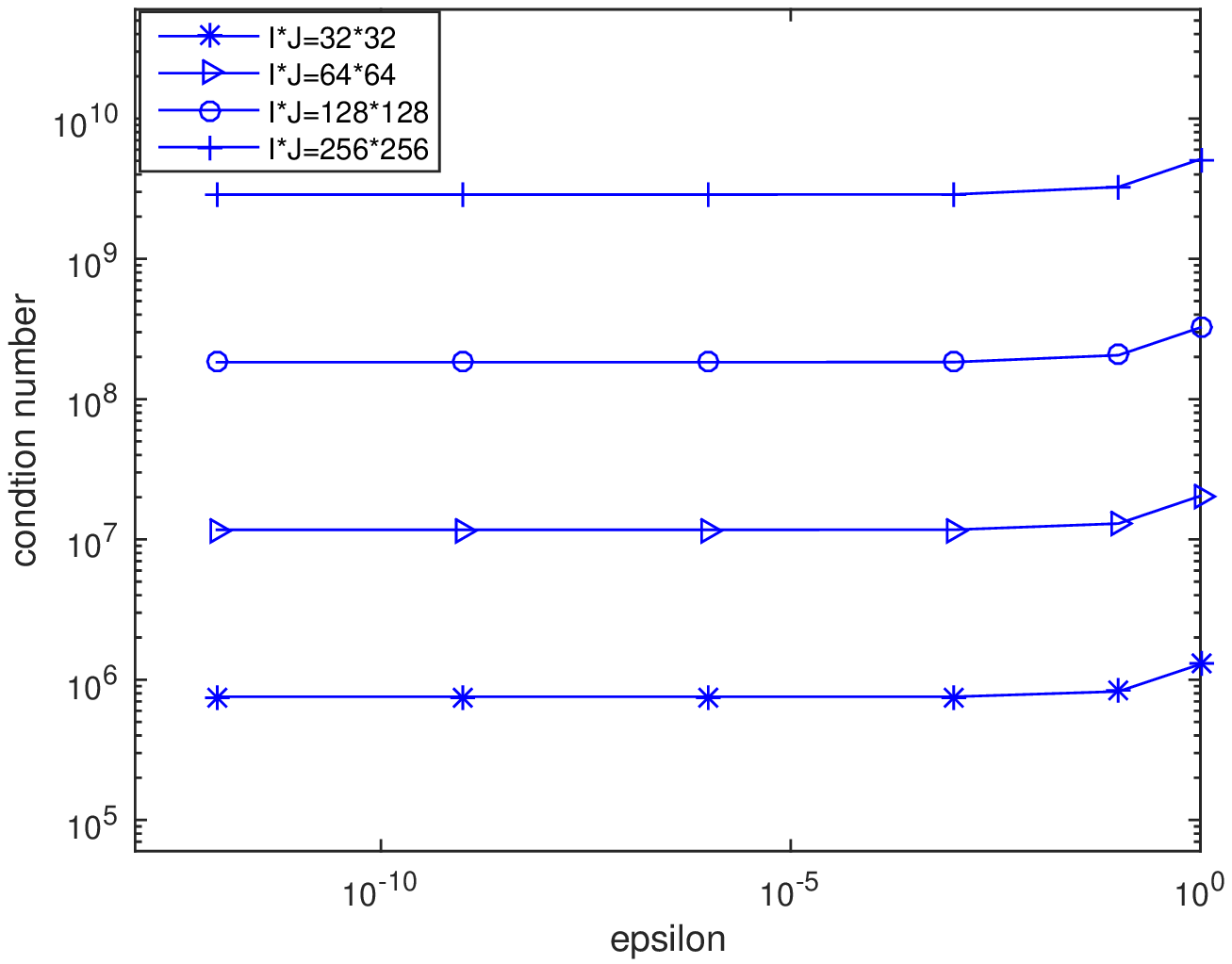} }
\subfigure [] {\includegraphics[width=7cm,clip]{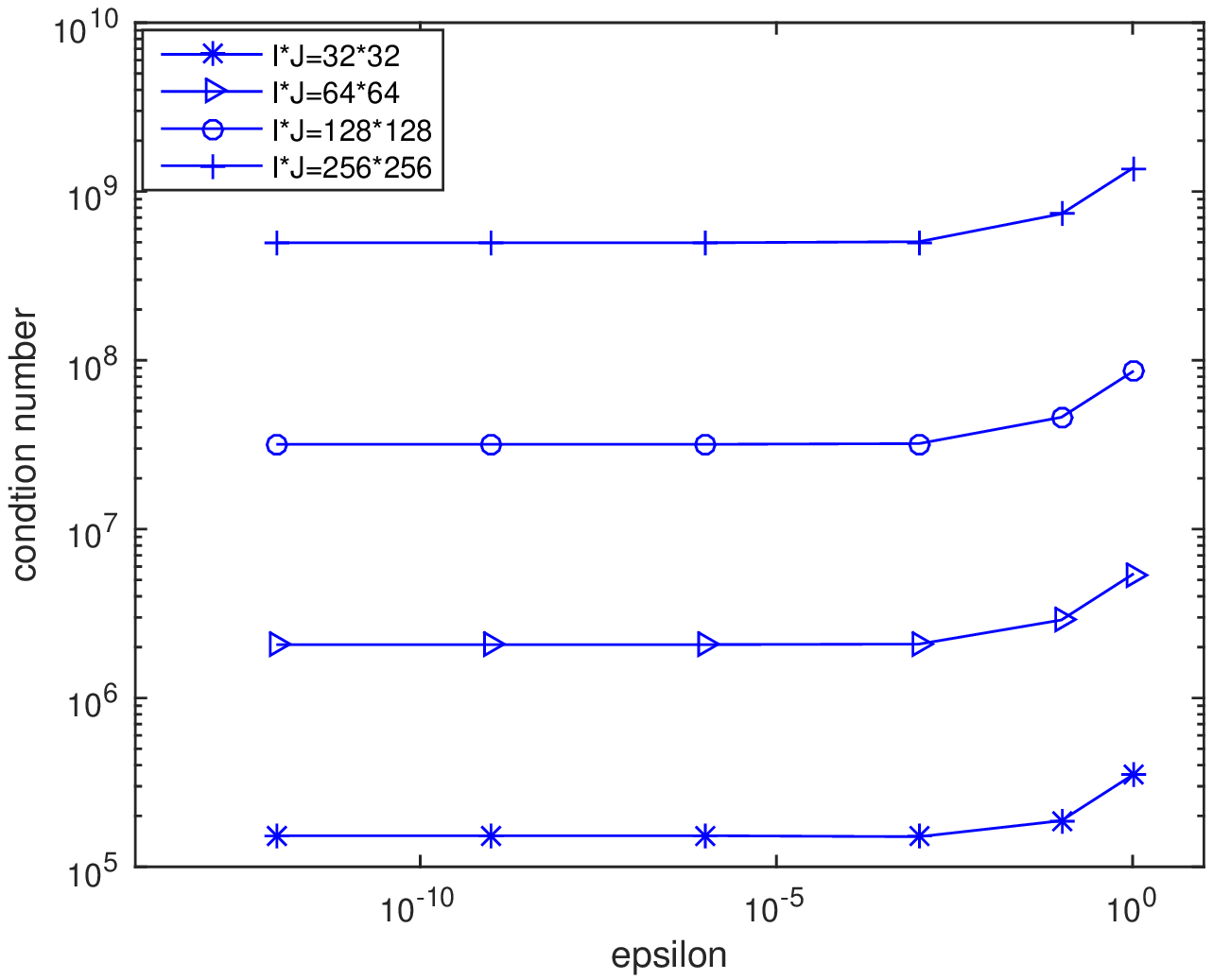} }
\subfigure[] {\includegraphics[width=7cm,clip]{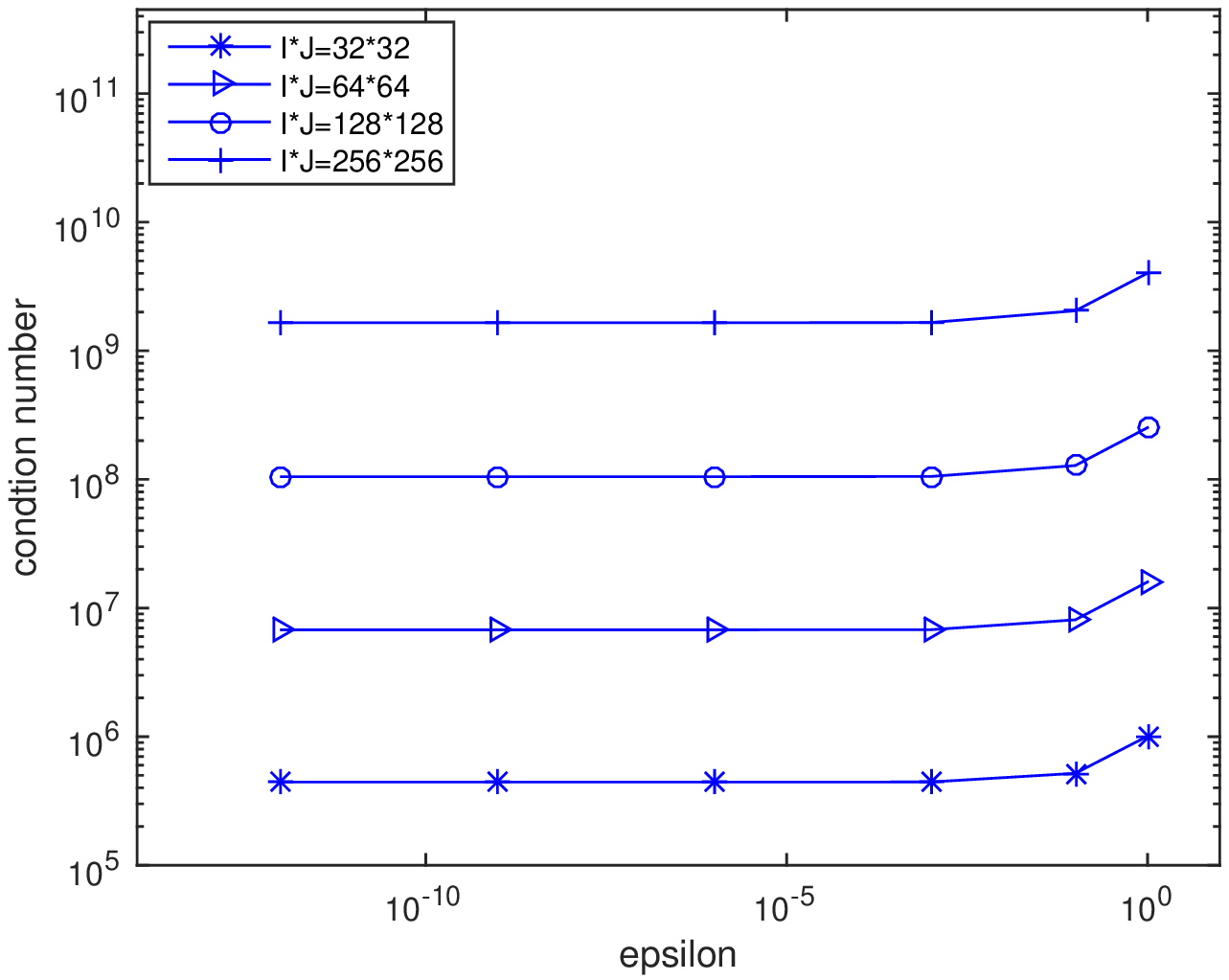}  }
\subfigure [] {\includegraphics[width=7cm,clip]{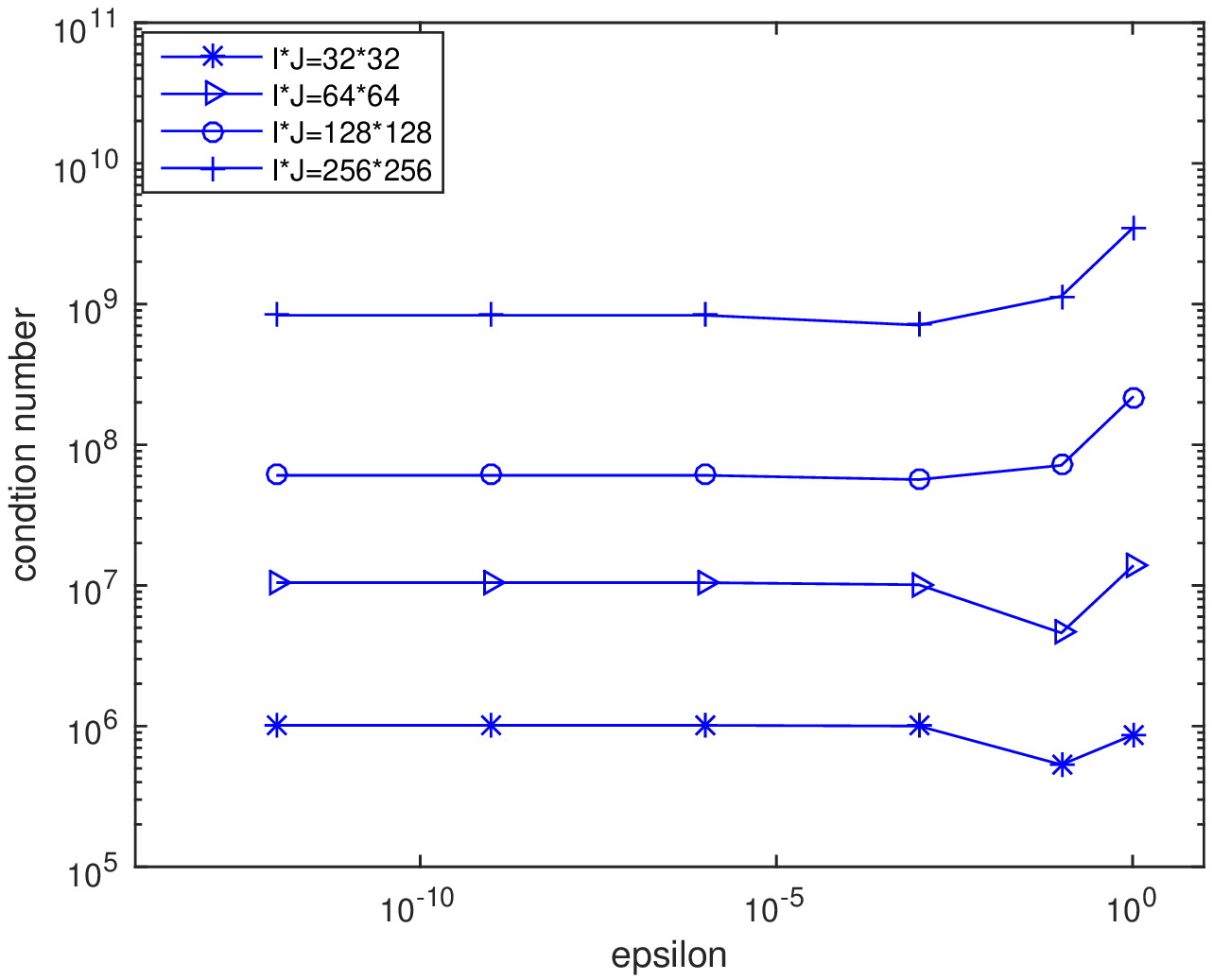} }

}
\caption{Example 1. The condition number for different grids and different $\epsilon$ values.
 (a). $\gamma_1=\gamma_2=0.5,~\varphi=0$.
 (a). $\gamma_1=0.5,~\gamma_2=0.85,~\varphi=0$.
 (b). $\gamma_1=0.5,~\gamma_2=0.85,~\varphi=\pi/4$.
 (c). $\gamma_1=0.25,~\gamma_2=0.85,~\varphi=\pi/3$.}\label{fig:ex1cond}
\end{figure}

\textbf{Condition Number }
The condition numbers of four cases are displayed on  Figure \ref{fig:ex1cond}, it is shown that the condition number of the discretized system is bounded by a constant independent of $\epsilon$ that has the same magnitude as the classical $9$-Point FDM, see Figure \ref{fig:ex1d1} $(b)$. However, \cite{Narski13} has pointed out that numerical discretizations of the original problem using magnetic field aligned coordinates lead to very badly conditioned systems when $\epsilon$ becomes small.

\begin{figure}[htb]
\centering
{
\subfigure [] {\includegraphics[width=7cm,clip]{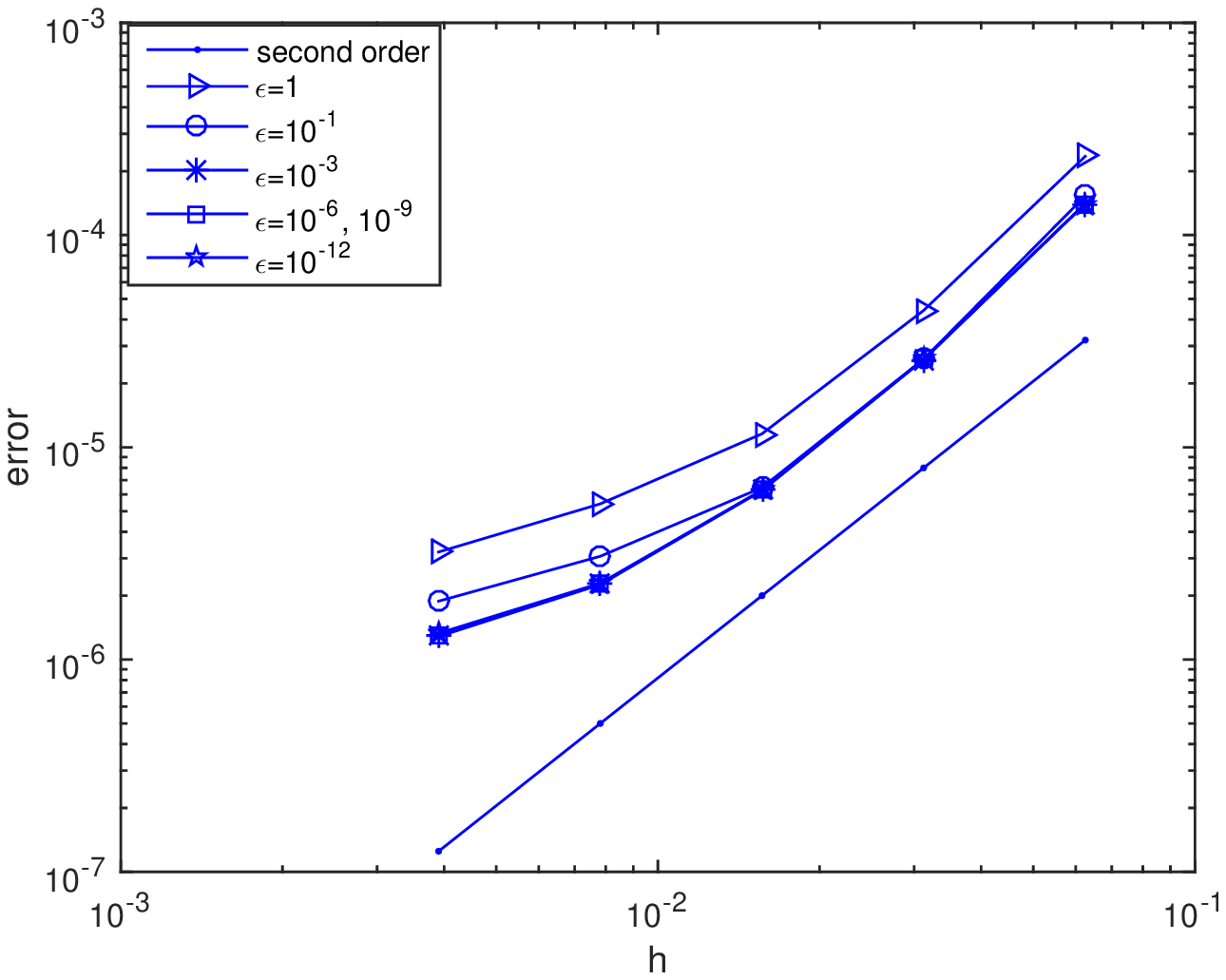} }
\subfigure[] {\includegraphics[width=7cm,clip]{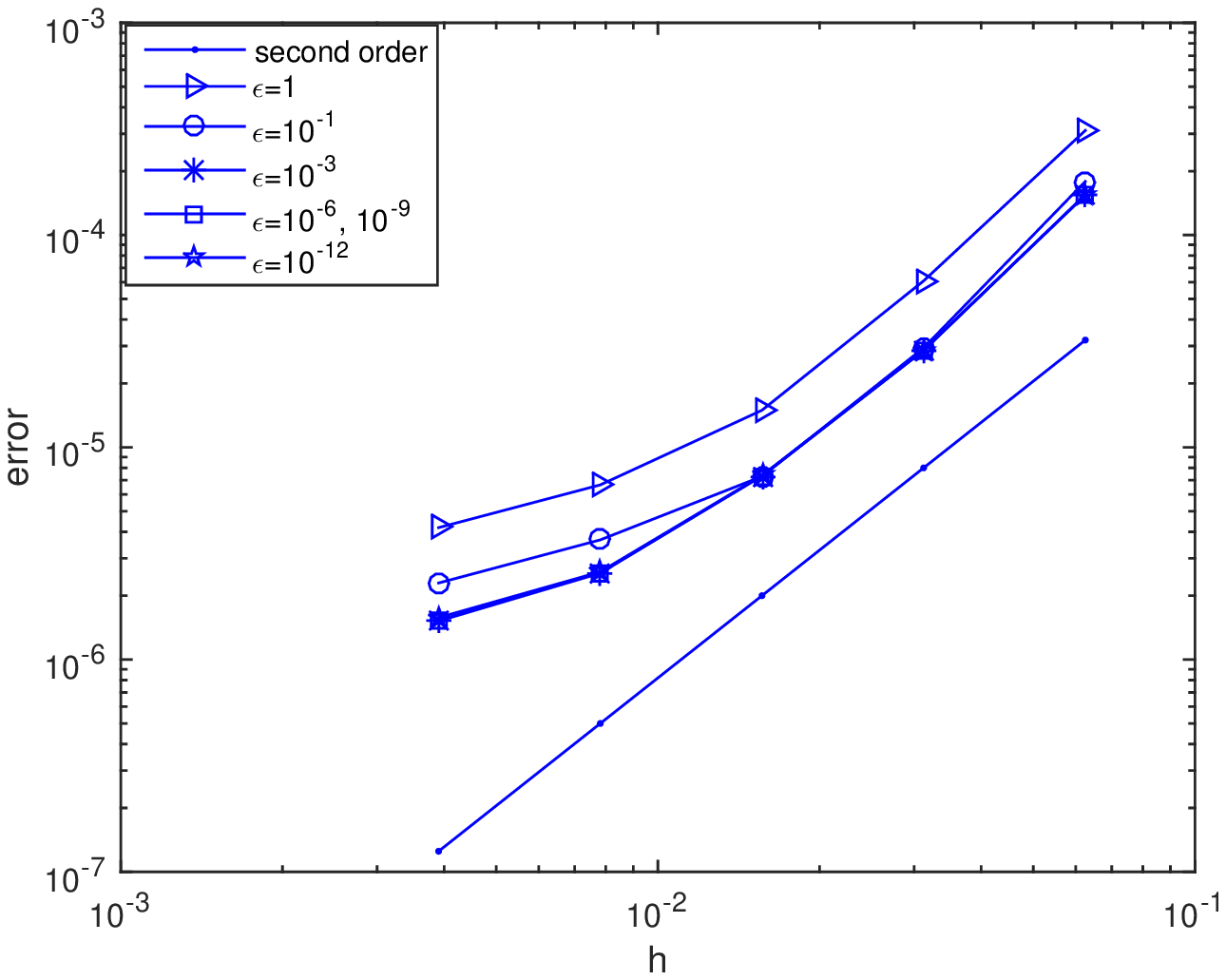}  }

}
\caption{Example 1. $\gamma_1=0.5,~\gamma_2=0.85,~\varphi=\pi/4$.  Convergence orders based on the closed field line is solved by "Method One". (a). Convergence orders  for $L^{2}$ error.
 (b). Convergence orders  for $L^{\infty}$ error.}\label{fig:ex1d2}
\end{figure}

As discussed in section 3.1, the closed field lines can be numerically determined in two different ways. We illustrate here that it is important to use "Method Two", in order to get uniform second order convergence.   Take $\gamma_1=0.5$, $\gamma_2=0.85$ and $\theta=\pi/4$ in \eqref{eq:uex} as an example. The results are displayed in Figure \ref{fig:ex1d2} and 
 no uniform convergence can be observed for "Method One", while "Method Two" works quite well.

\bigskip
\textbf{Example 2:}  We consider the case when there exhibit two ``magnetic islands''. Let the computational domain be $[-1,1]\times[-0.5,0.5]$ and exact solution be
 \beq\label{eq:ETwo}
 u^\epsilon_{exact}=\cos(\lambda \cos(2\pi(x-3/2))+\cos(\pi y))+\epsilon\sin(2\pi y)\sin(\pi x).
 \eeq
As $\epsilon \to 0$, $u^\epsilon_{excat}$ converges to $u^0=\cos(\lambda \cos(2\pi(x-3/2))+\cos(\pi y))$.
Let \beq\label{eq:test2B}
\mathbf{b}=\frac{B}{\mid B\mid},\quad B=\left(\begin{array}{c}-\pi\sin(\pi y)\\\lambda 2\pi\sin(2\pi(x-3/2))\end{array}\right)=\left(\begin{array}{c}B_1\\B_2\end{array}\right).
\eeq
Then $\mathbf{b}$ satisfies $\mathbf{b}\cdot\nabla u^0=0$, which indicates that the limiting solution $u^0$ is a constant along the field line. This is how we construct $\mathbf{b}$.

\begin{figure}[htb]
\centering
{
\subfigure [] {\includegraphics[width=7.5cm,clip]{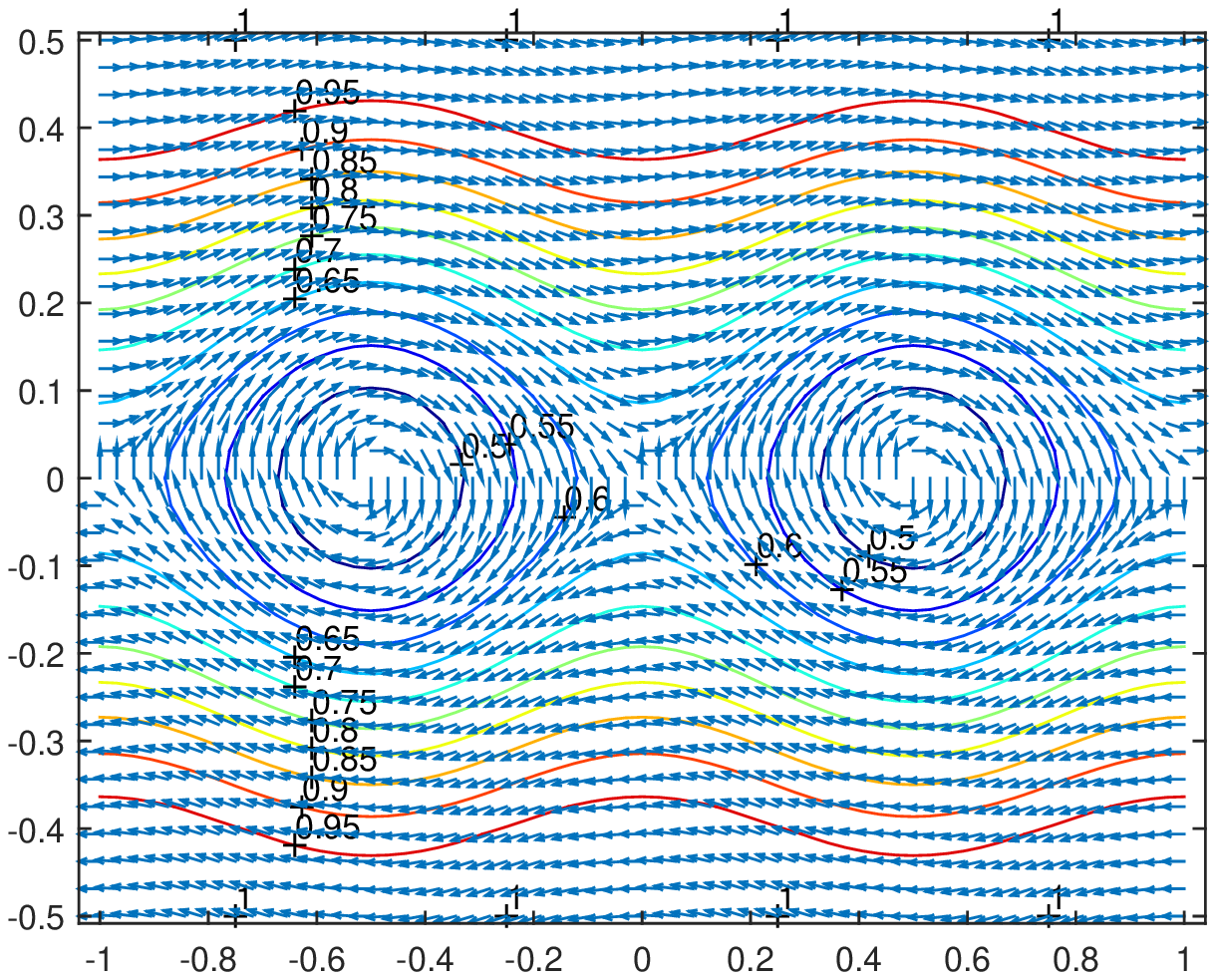} }
\subfigure[] {\includegraphics[width=7.5cm,clip]{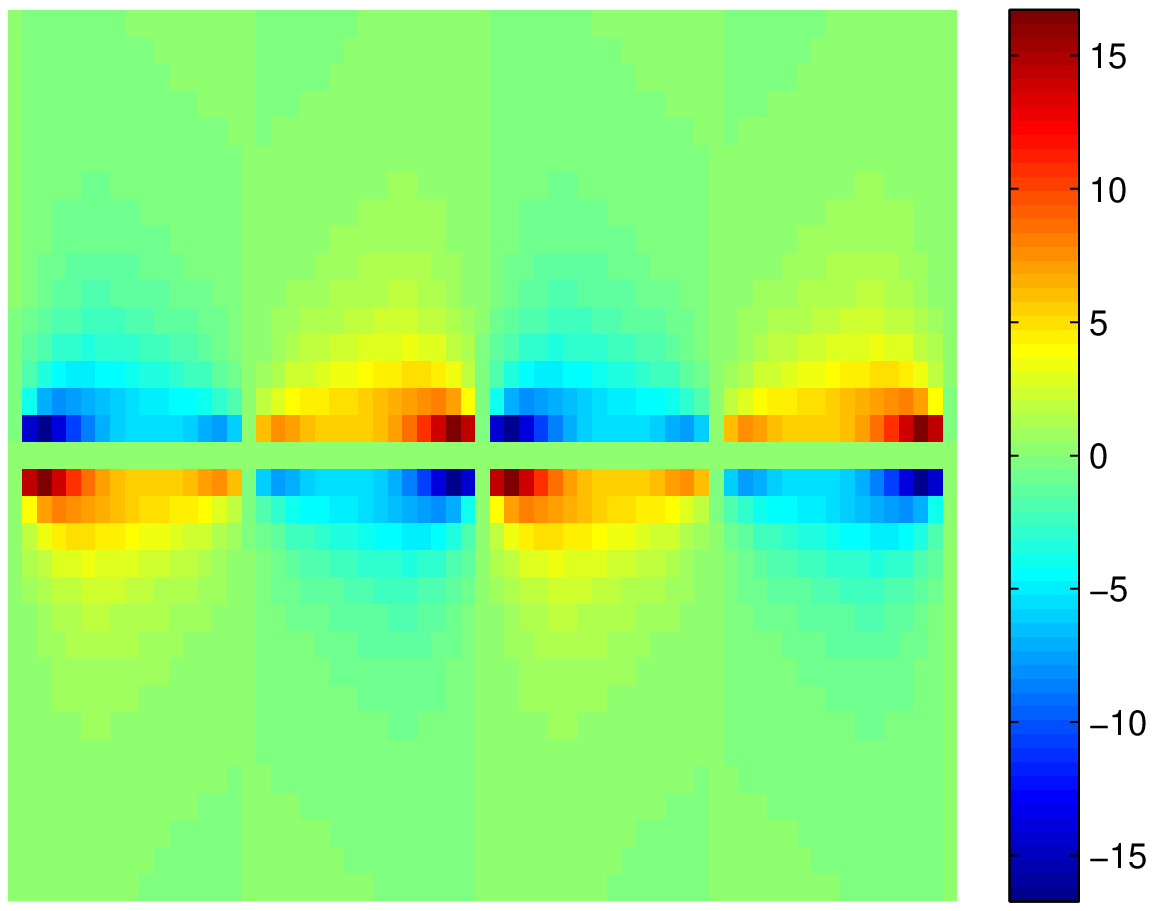}  }
\subfigure [] {\includegraphics[width=7.5cm,clip]{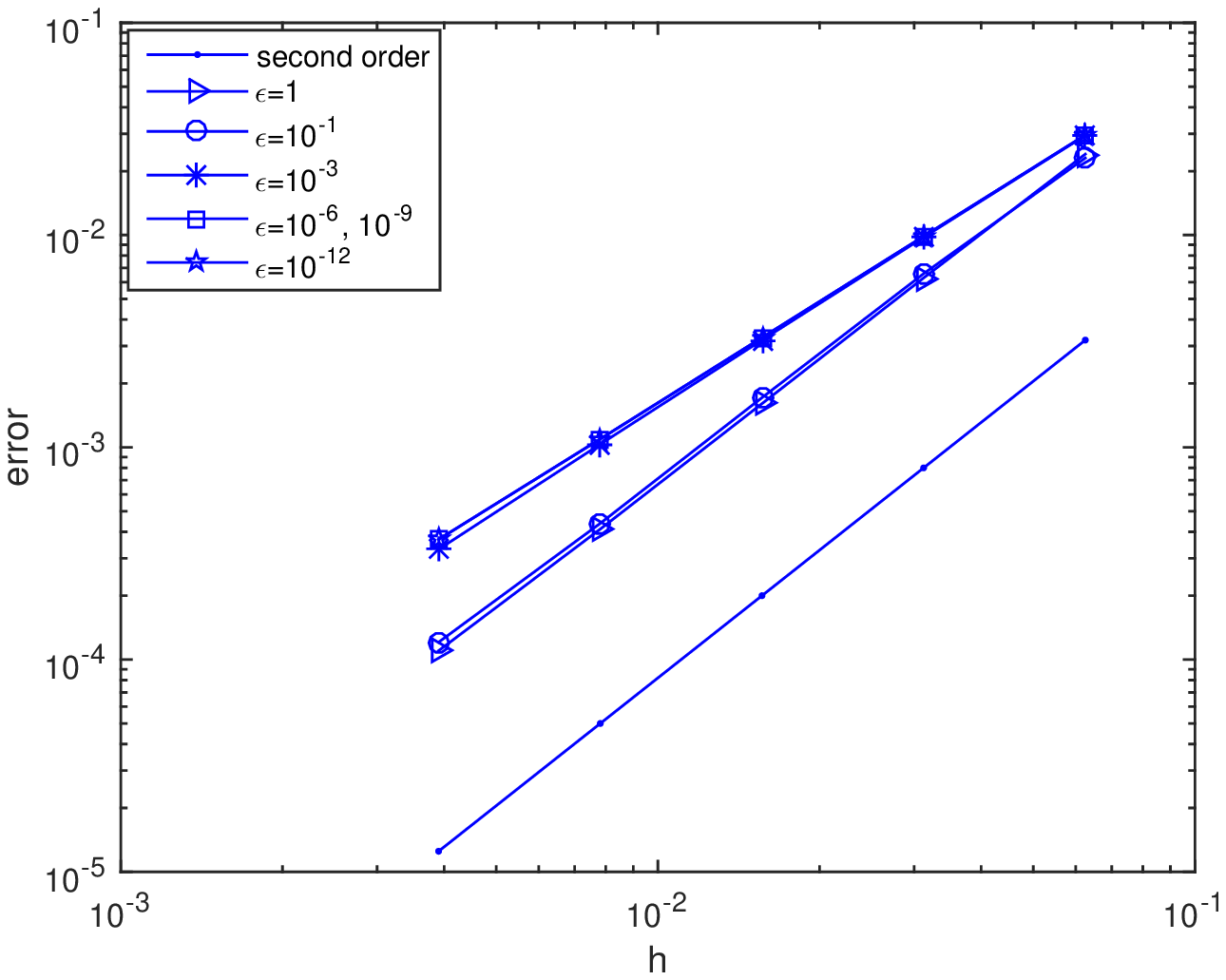} }
\subfigure[] {\includegraphics[width=7.5cm,clip]{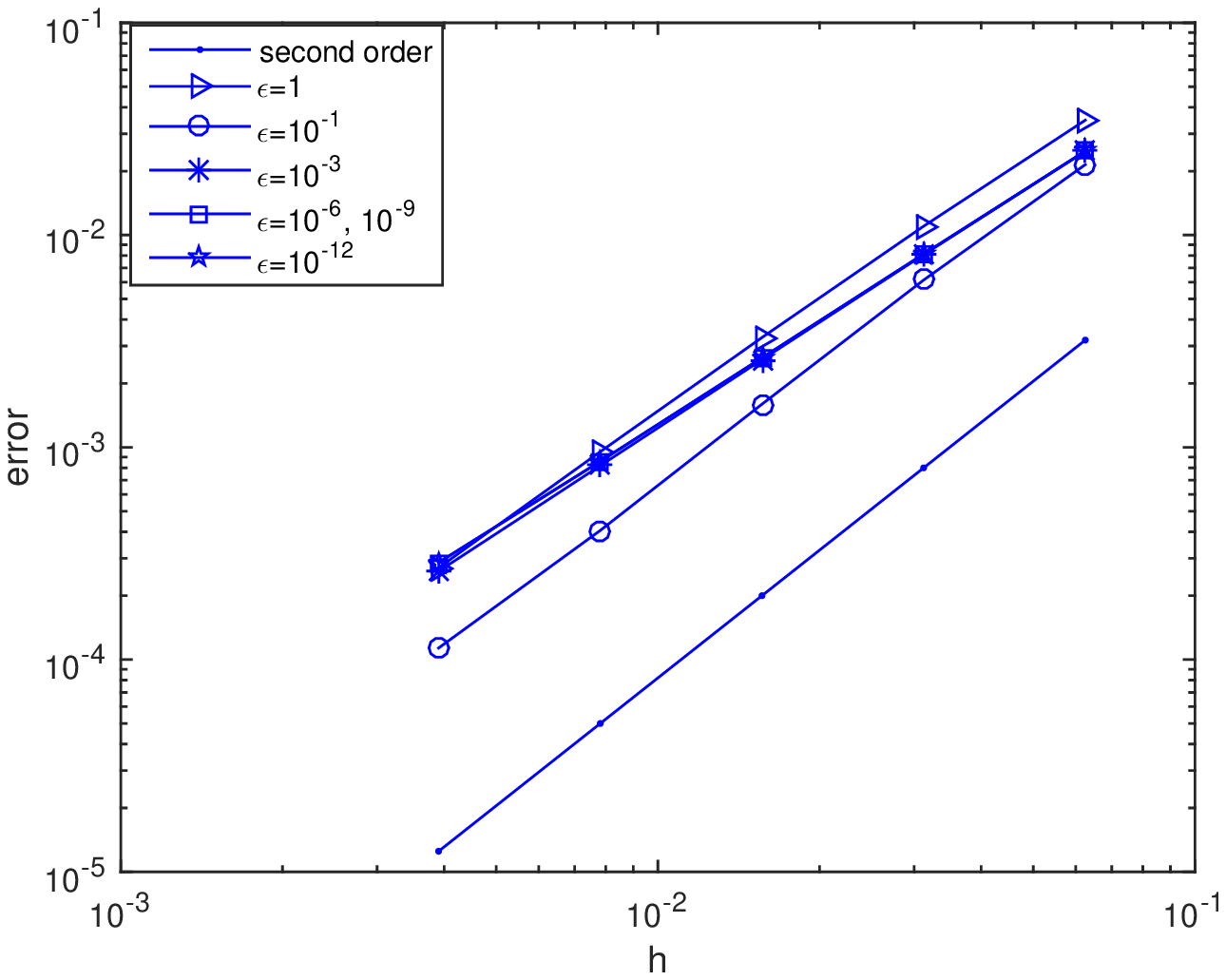}  }


}
\caption{ Example 2. $\lambda=0.1$.
 (a). Closed field.
(b). $\nabla\cdot\mathbf b$, $I\times J=64\times 32$.
 (c). Convergence orders  for $L^{2}$ error.
(d).  Convergence orders  for $L^{\infty}$ error.
 }
 \label{fig:ex2E}
\end{figure}

The source term is calculated by plugging \eqref{eq:ETwo}, \eqref{eq:test2B} into \eqref{eq:ellipT}. The field is showed in Figure \ref{fig:ex2E} $(a)$. In this example, $ \mathbf{b} $ has no definition at five points: $(-1,0)$, $(-0.5,0)$, $(0,0)$, $(0.5,0)$, $(1,0)$ and $\nabla \mathbf\cdot \mathbf b$ change quickly near these points (see Figure \ref{fig:ex2E} $(b)$).
Two of these five singular points $(-0.5,0)$ and $(0.5,0)$ are the centers of the ``magnetic island", while the other three are the connection points of two ``magnetic islands". ${\sqrt{B_1^2+B_2^2+\delta}}$  ($\delta=10^{-16}$) are used to replace ${\sqrt{B_1^2+B_2^2}}$, in order to approximate the values at the five singular points.
The convergence order of our scheme is around $1.5\sim1.6$ for $\epsilon$ ranging from $10^{-12}$ to $10^{-3}$ (See Figure \ref{fig:ex2E} $(c), (d)$). This is caused by the three singular points of $\mathbf{b}$ at $(-1,0)$, $(0,0)$, $(1,0)$.


%

\section{Conclusion}\label{con}

We present a simple Asymptotic-Preserving formulation for strongly anisotropic diffusion 
equation with closed field lines. The key idea is that we cut each of the closed field 
line at some point $(x_{0}, y_{0})$ and replace locally discretizing the equation at 
the point $(x_{0}, y_{0})$ by integration of the  differential equation along the cut
field line so that the singular $1/\epsilon$ terms disappear. The new system removes 
the ill-posedness and uniform second order convergence with respect to the anisotropy 
is observed numerically, even for the tilted ellipse case.

The scheme is efficient, general and easy to implement. The idea can be coupled with most standard discretizations and the computational cost keeps almost the same. Only slight modification to the original code is required, which makes it attractable to engineers.

\section*{Acknowledgement}
This work  was partially supported by NSFC 11301336
and 91330203. The authors would like to thank Professor Chunjing Xie for his valuable suggestions.

\end{document}